\documentclass{amsart}

\usepackage{amssymb, amscd}
\usepackage[all]{xy}
\usepackage{graphicx}
\usepackage{hyperref}

\numberwithin{equation}{section}

\newtheorem{thm}{Theorem}[section]
\newtheorem{prop}[thm]{Proposition}

\newtheorem{cor}[thm]{Corollary}
\newtheorem{conj}[thm]{Conjecture}
\newtheorem{q}[thm]{Question}

\theoremstyle{definition}

\theoremstyle{remark}
\newtheorem{rem}[thm]{Remark}
\newtheorem{exmp}[thm]{Example}

\renewcommand{\hom}{\operatorname{Hom}}
\renewcommand{\ker}{\operatorname{Ker}}

\newcommand{\Z}{\mathbb{Z}}
\newcommand{\Q}{\mathbb{Q}}
\newcommand{\R}{\mathbb{R}}
\newcommand{\C}{\mathbb{C}}
\newcommand{\F}{\mathbb{F}}

\DeclareMathOperator{\Ann}{Ann}

\DeclareMathOperator{\coker}{Coker}

\DeclareMathOperator{\GL}{GL}
\DeclareMathOperator{\im}{Im}
\DeclareMathOperator{\inj}{inj}

\DeclareMathOperator{\PSL}{PSL}
\DeclareMathOperator{\rank}{rank}

\DeclareMathOperator{\SL}{SL}
\DeclareMathOperator{\Spin}{Spin}

\DeclareMathOperator{\vol}{Vol}

\DeclareFontFamily{U}{mathx}{\hyphenchar\font45}
\DeclareFontShape{U}{mathx}{m}{n}{
      <5> <6> <7> <8> <9> <10>
      <10.95> <12> <14.4> <17.28> <20.74> <24.88>
      mathx10
      }{}
\DeclareSymbolFont{mathx}{U}{mathx}{m}{n}
\DeclareFontSubstitution{U}{mathx}{m}{n}
\DeclareMathAccent{\widecheck}{0}{mathx}{"71}

\begin{document}

%%%%%%% Title %%%%%%%%%%%%%%%%%%%%%%%%%%%%%%%%%%%%%%%%%%%%%%%%%%%%%%%%%%%%%%%%%
\title[A survey of the Thurston norm]{A survey of the Thurston norm}
\author[T.~Kitayama]{Takahiro Kitayama}
\address{Graduate School of Mathematical Sciences, the University of Tokyo, Japan}
\email{kitayama@ms.u-tokyo.ac.jp}

\begin{abstract}
We present an overview of the study of the Thurston norm, introduced by W.~P.~Thurston in the seminal paper ``A norm for the homology of $3$-manifolds'' (written in 1976 and published in 1986).
We first review fundamental properties of the Thurston norm of a $3$-manifold, including a construction of codimension-$1$ taut foliations from norm-minimizing embedded surfaces, established by D.~Gabai.
In the main part we describe relationships between the Thurston norm and other topological invariants of a $3$-manifold:
the Alexander polynomial and its various generalizations, Reidemeister torsion, the Seiberg--Witten invariant, Heegaard Floer homology, the complexity of triangulations and the profinite completion of the fundamental group.
Some conjectures and questions on related topics are also collected.

The final version of this paper will appear as a chapter in the book ``In the tradition of Thurston, II'', edited by K.~Ohshika and A.~Papadopoulos (Springer, 2022).
\end{abstract}

\maketitle

\setcounter{tocdepth}{1}
\tableofcontents

%%%%%%%% Section 1 %%%%%%%%%%%%%%%%%%%%%%%%%%%%%%%%%%%%%%%%%%%%%%%%%%%%%%%%%%%%
\section{Introduction} \label{sec:1}

In the seminal paper \cite{Th86} written in 1976, Thurston introduced a seminorm on the first real cohomology group of a $3$-manifold, called the \textit{Thurston norm}.
It measures the topological complexity of embedded surfaces dual to a given integral cohomology class.
The unit ball of the seminorm constitutes a convex polyhedron and Thurston described the distribution of the cohomology classes represented by fibrations over a circle in terms of the top-dimensional faces of the polyhedron.
Thurston also showed that every compact leaf of a codimension-$1$ taut foliation minimizes the Thurston norm.
The Thurston norm has become a fundamental tool in the study of incompressible surfaces, fibrations over a circle and codimension-$1$ foliations.

The roots of the Thurston norm go back to the study of the genus of a knot in the $3$-sphere, and it has been extensively studied in various contexts up to the present.
The goal of this chapter is to present an overview of the study of the Thurston norm, with an emphasis on its relationships with other topological invariants, and without aiming for completeness.
We summarize some of the results discussed in this chapter in the following.

Inspired by the foundational work~\cite{Th86} of Thurston, Gabai~\cite{Ga83} developed sutured manifold theory, extended by Scharlemann~\cite{Sc89}, and established a construction of codimension-$1$ taut foliations from embedded surfaces minimizing the Thurston norm. 
As consequences Gabai resolved the Property R conjecture and the Po\'enaru conjecture, and showed the equivalence of the Thurston norm and the Gromov norm~\cite{Gro82} on the second homology group.
Sutured manifold theory provides an efficient algorithm to compute the Thurston norm.
Tollefson and Wang~\cite{TW94, TW96}, and later Cooper and Tillmann~\cite{CT09} described another algorithm via normal surface theory.

One of the most fundamental algebraic invariants related to the Thurston norm is the \textit{Alexander polynomial}, equivalent to \textit{Milnor torsion}~\cite{Mi66, Tu75, Tu02b}.
It is well known that the classical Alexander polynomial of a knot in the $3$-sphere gives a lower bound on its genus.
As a generalization, McMullen~\cite{Mc02} introduced the \textit{Alexander norm} for a general $3$-manifold, and showed that it gives a lower bound on the Thurston norm.
Also, gauge theory has been a successful tool to study the complexity of embedded surfaces in $3$- and $4$-manifolds.
In the \textit{adjunction inequality} for a $3$-manifold the \textit{Seiberg--Witten invariant}~\cite{Wit94} gives a lower bound on the Thurston norm~\cite{Au96, Kr98}.
By the equivalence of Milnor torsion and the Seiberg--Witten invariant of a $3$-manifold~\cite{MeTa96, Tu98}, the above two lower bounds coincide~\cite{Kr98, V03}.
Furthermore, Kronheimer and Mrowka~\cite{KrMr97b} described the Thurston norm in terms of solutions of the Seiberg--Witten monopole equations.

These relationships were generalized in \textit{Heegaard Floer homology}~\cite{OS04b, OS04e} and \textit{monopole Floer homology}~\cite{KrMr07}, which provide categorifications of the Milnor torsion and the Seiberg--Witten invariant of a $3$-manifold.
These homology theories were also shown to be equivalent~\cite{CGH11, CGH12a, CGH12b, KLT20a, KLT20b, KLT20c, KLT20d}.
Ozsv\'ath and Szab\'o~\cite{OS04b} showed that Heegaard Floer homology determines the Thurston norm of a closed $3$-manifold, and Ni~\cite{Ni09b} showed that Heegaard Floer homology detects fiberedness of a closed $3$-manifold.
Knot Floer homology~\cite{OS04c, Ra03} provides a categorification of the classcal Alexander polynomial of a knot.
Ozsv\'ath and Szab\'o~\cite{OS04b} showed that knot Floer homology determines the genus of a knot, and Ghiggini~\cite{Gh08}, Ni~\cite{Ni07}, and Juh\'asz~\cite{Ju08, Ju10} showed that knot Floer homology detects fiberedness of a knot.

Twisted Alexander polynomials~\cite{Linx01, Wa94} associated with linear representations, and higher-order Alexander polynomials~\cite{Co04, Ha05} with coefficients in skew-fields are more direct generalizations of the Alexander polynomial.
These polynomials give generalized lower bounds on the Thurston norm~\cite{Fril07, FrKim08b, HF07}.
Furthermore, Friedl and Vidussi~\cite{FV08, FV11a, FV11c, FV14a} showed that twisted Alexander polynomials detect fiberedness of a $3$-manifold, and Friedl, Nagel and Vidussi~\cite{FN15, FV15} showed that twisted Alexander polynomials determine the Thurston norm.
The $L^2$-Alexander invariant or torsion~\cite{DFL16, LZ06} are ``polynomial-like'' $L^2$-invariants, generalizing the $L^2$-torsion~\cite{Lu02}.
Friedl and L\"uck~\cite{FL19b}, and Liu~\cite{Liu17} showed that the $L^2$-Alexander torsion determines the Thurston norm.

Boileau and Friedl~\cite{BoFr20}, Bridson, Reid and Wilton~\cite{BR20, BRW17}, and Liu~\cite{Liu20} showed certain rigidity results of the Thurston norm and fiberedness of a $3$-manifold on the profinite completion of the fundamental group.
Also, Jaco, Rubinstein, Spreer and Tillmann~\cite{JRT13, JRST20a, JRST20b} introduced a $\Z / 2 \Z$-analogue of the Thurston norm and showed that it gives lower bounds on minimal numbers of tetrahedra in triangulations and ideal triangulations of a $3$-manifold. 

The influence of the Thurston norm is not limited to low-dimensional topology, and the following significant topics are, for example, unfortunately beyond the scope of this article.
The universal $L^2$-torsion defines an equivalence class of a pair of convex polytopes for a torsion-free group satisfying certain conditions.
Such an equivalence class can be regarded as the unit ball of the (dual) Thurston norm, and as already shown in~\cite{GKL21, HK18, HK20, FL17, FLT19, FL19a, FST17, FT20, FuKie18, Kie20a}, there should be a fruitful theory for the ``Thurston norm of groups''.
Calegari~\cite{Ca08, Ca09a, Ca09b, Ca13, CaGo13} studied a group-theoretical interpretation of the Thurston norm in terms of the stable commutative length.
Also, Flores, Kahrobaei and Koberda~\cite{FlKaKo20} proposed a public-key and a symmetric-key cryptographic schemes based on the Thurston norm of hyperbolic $3$-manifolds.

For foundational results on the Thurston norm there are already excellent expositions in \cite{Ga98, Oe86, Sc89}, and also in \cite[Chapter 10]{CandC03} and \cite[Chapter 2]{Kap01}.
See also the survey \cite[Section 12]{Sa20} of the impact of Thurston's work on knot theory in the first volume of this series of books.
For terminology and developments of the study of $3$-manifolds we refer the reader to the book~\cite{AFW15}.

Throughout we do not attempt to state results in their greatest generality, and we do not make any claims to originality.

\subsection*{Organization}

Section~\ref{sec:2} provides a brief review of the definition and fundamental properties of the Thurston norm, including the correspondence between embedded surfaces minimizing the Thurston norm and codimension-$1$ taut foliations. 
Section~\ref{sec:3} describes the relationships between the norms on the first cohomology group associated with the Alexander and Teichm\"uller polynomials and the Thurston norm.
Section~\ref{sec:4} discusses adjunction inequalities from Seiberg--Witten theory for $3$- and $4$-manifolds.
Section~\ref{sec:5} summarizes the facts that Heegaard Floer homology and knot Floer homology detect the Thurston norm, knot genus and fiberedness of a $3$-manifold and a knot.
Section~\ref{sec:6} deals with twisted Alexander polynomials, higher-order Alexander polynomials and $L^2$-Alexander torsion.
Here results on the Thurston norm are described in terms of Reidemeister torsion.
Section~\ref{sec:7} contains constructions of the Thurston norm ball via normal surface theory and applications of a $\Z / 2 \Z$-analogue of the Thurston norm to the study of complexity of a $3$-manifold.
Section~\ref{sec:8} is devoted to explain certain rigidity results of the Thurston norm on the profinite completion of the fundamental group.
In Section~\ref{sec:9} we conclude by collecting some conjectures and questions on the Thurston norm and related topics.

\subsection*{Conventions and notation}

All surfaces and manifolds are understood to be compact, connected and oriented unless we say specifically otherwise.
For a link $L$ in $S^3$ we denote by $X_L$ the complement of an open tubular neighborhood of $L$.
For an integral domain $R$ we denote by $Q(R)$ its quotient field.

\subsection*{Acknowledgments}

The author wishes to thank Ken'ichi Ohshika and Athanase Papadopoulos for inviting him to write this survey article and for a careful check of it greatly improving the presentation.
The author also wishes to thank Stefan Friedl for many helpful suggestions, and Martin Bridson, Alan Reid and Makoto Sakuma for valuable comments.
The author was supported by JSPS KAKENHI Grant Numbers JP18K13404, JP18KK0071, JP18KK0380 and JP21H00986.

%%%%%%%% Section 2 %%%%%%%%%%%%%%%%%%%%%%%%%%%%%%%%%%%%%%%%%%%%%%%%%%%%%%%%%%%%
\section{Foundations of the Thurston norm} \label{sec:2}

First we briefly review the definition of the Thurston norm of a $3$-manifold and its fundamental properties.
We summarize original results by Thurston~\cite{Th86} and Gabai~\cite{Ga83, Ga87} on the polyhedron structure of the unit ball of the Thurston norm, the distribution of cohomology classes represented by fibrations over a circle, the correspondence between embedded surfaces minimizing the Thurston norm and codimension-$1$ taut foliations, and the equivalence of the Thurston and Gromov norms.

\subsection{Thurston norm}

We begin with the definition of the Thurston norm of a $3$-manifold $M$~\cite{Th86}.

For a surface $S$ with connected components $S_1$, $S_2$, $\dots$, $S_k$ its \textit{complexity} $\chi_-(S)$ is defined by
\[ \chi_-(S) = \sum_{i = 1}^k \max \{ -\chi(S_i), 0 \}, \]
where $\chi$ is the Euler characteristic.
Every cohomology class in $H^1(M; \Z)$ is represented by a smooth map $M \to S^1$ and the properly embedded surface obtained as the inverse image of any regular value represents the Poincar\'e dual of the cohomology class.
The \textit{Thurston norm} $x_M$ on $H^1(M; \Z)$ of $M$ is defined by 
\[ x_M(\phi) = \min \{ \chi_-(S) ~;~ \text{$S$ is a properly embedded surface in $M$ dual to $\phi$} \} \]
for $\phi \in H^1(M; \Z)$.

In \cite{Th86} Thurston first showed that $x_M$ is a seminorm on $H^1(M; \Z)$.
Key observations are that the $k$-multiple of a homology class is represented by $|k|$ disjoint properly embedded surfaces representing the homology class, and that the ``double curved sum'' of two properly embedded surfaces with transverse intersection represents the sum of their homology classes.
Since $x_M$ is linear on each ray through the origin, it extends to $H^1(M; \Q)$.
Since $x_M$ is a convex function, it extends to all of $H^1(M; \R)$ in a unique continuous way.
We denote also by $x_M$ the extended seminorm on $H^1(M; \R)$.
Moreover, for $\phi \in H^1(M; \R)$ with $x_M(\phi) = 0$ the ray through $\phi$ comes arbitrarily near lattice points, and if nonzero multiple $a \phi$ is near enough to a lattice point $l$, then the integer $x_M(l) = x_M(l - a \phi)$ must be $0$.
Thus $\phi$ can be approximated by multiples of lattice points $l$ with $x_M(l) = 0$.
Summarizing, we state the following theorem~\cite[Theorem 1]{Th86}:

\begin{thm}[\cite{Th86}] \label{thm:norm}
The Thurston norm $x_M$ uniquely extends to $H^1(M; \R)$ as a seminorm:
\begin{enumerate}
\item $x_M(a \phi) = |a| x_M(\phi)$,
\item $x_M(\phi + \psi) \leq x_M(\phi) + x_M(\psi)$,
\end{enumerate}
for $\phi$, $\psi \in H^1(M; \R)$ and $a \in \R$.
Moreover, $x_M^{-1}(\{ 0 \})$ is spanned by integral cohomology classes dual to properly embedded surfaces in $M$ with non-negative Euler characteristic.
\end{thm}

\begin{rem}
More generally, the seminorm can be defined on $H_2(M, A; \R)$ for any submanifold $A$ in $\partial M$, as Scharlemann described in \cite{Sc89}.
\end{rem}

In general, $x_M$ is only a seminorm, but Theorem~\ref{thm:norm}, in particular, shows that  $x_M$ is nondegenerate for a \textit{hyperbolic $3$-manifold}, i.e., a $3$-manifold whose interior admits a complete Riemannian metric of constant sectional curvature $-1$ and finite volume.

A properly embedded surface $S$ is called \textit{norm-minimizing} if $\chi_-(S) = x_M(\phi)$ for its dual $\phi \in H^1(M; \Z)$.
Every connected norm-minimizing surface $S$ with negative Euler characteristic is incompressible since any compression of such a surface $S$ along a simple closed curve not bounding any disc in $S$ would reduce $\chi_-(S)$.

\begin{exmp}
The Thurston norm is a generalization of the knot genus:
The \textit{genus} $g(K)$ of a knot $K$ in $S^3$ is the minimum genus of Seifert surfaces of $K$.
Every norm-minimizing surface in the complement $X_K$ of $K$ dual to a generator $\psi \in H^1(X_K; \Z)$ corresponds to a minimal genus Seifert surface of $K$, and we have
\[ x_{X_K}(\phi) = 2 g(K) - 1 \]
for a nontrivial knot $K$.
\end{exmp}

\begin{exmp}
Let $M$ be a $3$-manifold fibering over a circle with a fiber surface $S$.
Then every incompressible surface in $M$ representing the same homology class in $H_2(M, \partial M; \Z)$ as $S$ is isotopic to $S$.

We give a sketch of the proof as in \cite[Lemma 5.1]{EL83}.
First such an incompressible surface $S'$ in $M$ lifts homeomorphically to an incompressible surface $\overline{S}'$ in the infinite cyclic covering $S \times \R$ of $M$ corresponding with the fibration.
Then the inclusion-induced homomorphism $\pi_1 \overline{S}' \to \pi_1 (S \times \R)$ is an isomorphism.
If it would be not surjective, then an argument with van Kampen's theorem would imply that $\pi_1 (S \times \R)$ is not finitely generated.
We thus see that $S'$ is isotopic to $S$.

An immediate consequence is that $S$ is norm-minimizing.
We will see a more general result in Theorem~\ref{thm:foliation1} for codimension-$1$ foliations on $M$.
\end{exmp}

\subsection{Norm balls and fibrations over a circle}

We next discuss the structure of the unit ball of the Thurston norm and the distribution of the cohomology classes represented by fibrations over a circle.

The \textit{Thurston norm ball} of a $3$-manifold $M$, denoted by $B_M$, is the unit ball of $x_M$:
\[ B_M = \{ \phi \in H^1(M; \R) ~;~ x_M(\phi) \leq 1 \}. \]
A seminorm determines its unit ball and vice versa.
We set
\[ \Ann(x_M^{-1}(\{ 0 \})) = \{ \alpha \in H_1(M; \R) ~;~ \langle \phi, \alpha \rangle = 0 ~\text{for all $\phi \in x_M^{-1}(\{ 0 \}) $} \}, \]
where $\langle \cdot, \cdot \rangle$ is the Kronecker pairing.
Note that if $x_M$ is nondegenerate, then $\Ann(x_M^{-1}(\{ 0 \})) = H_1(M; \R)$.
The \textit{dual Thurston norm} $x_M^*$ on $\Ann(x_M^{-1}(\{ 0 \}))$ is defined by
\[ x_M^*(\alpha) = \sup \{ \langle \phi, \alpha \rangle ~;~ \phi \in B_M \} \]
for $\alpha \in \Ann(x_M^{-1}(\{ 0 \}))$.
Theorem~\ref{thm:norm} implies that $x_M^*$ is a norm on $\Ann(x_M^{-1}(\{ 0 \})$.
The \textit{dual Thurston norm ball} of $M$, denoted by $B_M^*$, is the unit ball of $x_M^*$:
\begin{align*} 
B_M^* &= \{ \alpha \in \Ann(x_M^{-1}(\{ 0 \}) ~;~ x_M^*(\alpha) \leq 1 \} \\
&= \{ \alpha \in \Ann(x_M^{-1}(\{ 0 \}) ~;~ \langle \phi, \alpha \rangle \leq 1 ~\text{for all $\phi \in B_M$} \}.
\end{align*}

The unit ball of a seminorm, a priori, may be an arbitrary convex body symmetric in origin,
but Thurston~\cite[Theorem 2]{Th86} showed that the structure of $B_M$ is more restrictive.
 
\begin{thm}[\cite{Th86}] \label{thm:polytope}
The dual Thurston norm ball $B_M^*$ of a $3$-manifold $M$ is a convex polytope in $H_1(M; \R)$ with finitely many vertices $\pm \alpha_1$, $\dots$, $\pm \alpha_k \in \Ann(x_M^{-1}(\{ 0 \})) \cap H_1(M; \Z)$, and we have
\[ B_M = \{ \phi \in H^1(M; \R) ~;~ |\langle \phi, \alpha_i \rangle| \leq 1 ~\text{for $1 \leq i \leq k$} \}. \] 
\end{thm}

Theorem~\ref{thm:polytope} is a formal consequence of the fact that $x_M$ is $\Z$-valued on the integral lattice $H^1(M; \Z)$.

\begin{cor}[\cite{Th86}]
The Thurston norm ball $B_M$ of a $3$-manifold $M$ is a (possibly noncompact) convex polyhedron in $H^1(M; \R)$ with finitely many vertices in $H^1(M; \Q)$.
\end{cor}

\begin{rem}
In this chapter a \textit{convex polyhedron} in a real affine linear space refers to a closed convex subset such that every point on the boundary lies in only finitely many maximal convex subsets of the boundary.
A \textit{convex polytope} refers to a compact convex polyhedron.
\end{rem}

A cohomology class $\phi \in H^1(M; \Z)$ is called \textit{fibered} if $M$ fibers over a circle such that the fibers are dual to $\phi$.
Since integration of a nonsingular closed $1$-form on $M$ with integer periods defines a fibration over a circle, $\phi \in H^1(M; \Z)$ is fibered if and only if $\phi$ is represented by a nonsingular closed $1$-form on $M$.

An observation is that since every nonsingular closed $1$-form on $M$ remains nonsingular after sufficiently small perturbation, the subset of cohomology classes of nonsingular closed $1$-forms is open in $H^1(M; \R)$.
Also, a nonsingular closed $1$-form on $M$ defines a codimension-$1$ foliation on $M$, which we will discuss in Section~\ref{ss:foliation}.
Based on the study of general position of incompressible surfaces with respect to codimension-$1$ foliations~\cite[Theorem 4]{Th86}, Thurston~\cite[Theorem 5]{Th86} described the distribution of fibered classes in terms of $B_M$ as follows:

\begin{thm}[\cite{Th86}] \label{thm:fibered}
Let $M$ be a $3$-manifold fibering over a circle with fiber of negative Euler characteristic.
There are some top-dimensional faces of $B_M$ such that $\phi \in H^1(M; \Z)$ is fibered if and only if $\phi$ lies in the interior of the cone on one of the faces.
\end{thm}

\begin{rem}
For a $3$-manifold $M$ fibering over a circle with fiber of nonnegative Euler characteristic, $x_M$ vanishes on $H^1(M; \R)$.
\end{rem}

Such top-dimensional faces of $B_M$ as in Theorem~\ref{thm:fibered} are called \textit{fibered faces} of $B_M$.

A $3$-manifold $M$ is \textit{atoroidal} if $M$ contains no incompressible torus.

\begin{cor}[\cite{Th86}] \label{cor:fibered}
Let $M$ be an atoroidal $3$-manifold with $b_1(M) > 1$.
Then there exists an incompressible surface which is not the fiber of a fibration over a circle.
\end{cor}

\begin{rem}
In the proof of Thurston's hyperbolization theorem~\cite{Ot96, Ot98, Kap01} Corollary~\ref{cor:fibered} played a significant role to reduce the exceptional (semi)fibered case to the generic case.
\end{rem}

For a $3$-manifold $M$ with $b_1(M) = 1$ its norm balls $B_M$ and $B_M^*$ are closed intervals centered at origins, possibly consisting only of origins.
The following examples together with one for the complement of the $3$-link chain in $S^3$ are given in \cite[Examples 1, 2, 3]{Th86}.
See also \cite[Section 4]{Th86} for a large variety of shapes for (dual) Thurston norm balls.

\begin{exmp}
Let $L$ be the Whitehead link.
Let $\mu_1$, $\mu_2 \in H_1(X_L; \R)$ be a basis represented by meridians of the two components of $L$, and $\mu_1^*$, $\mu_2^* \in H^1(X_L; \R)$ its dual basis.
Then $B_{X_L}$ is the diamond with vertices $\pm \mu_1^*$, $\pm \mu_2^*$ and $B_{X_L}^*$ is the square with vertices $\pm \mu_1 \pm \mu_2$.
All the $2$-dimensional faces of $B_{X_L}$ are fibered faces.
\end{exmp}

\begin{exmp}
Let $L$ be the Borromean rings.
Let $\mu_1$, $\mu_2$, $\mu_3 \in H_1(X_L; \R)$ be a basis represented by meridians of the three components of $L$, and $\mu_1^*$, $\mu_2^*$, $\mu_3^* \in H^1(X_L; \R)$ its dual basis.
Then $B_{X_L}$ is the octahedron with vertices $\pm \mu_1^*$, $\pm \mu_2^*$, $\pm \mu_3^*$ and $B_{X_L}^*$ is the cube with vertices $\pm \mu_1 \pm \mu_2 \pm \mu_3$.
All the $3$-dimensional faces of $B_{X_L}$ are fibered faces.
\end{exmp}

\subsection{Norm-minimizing surfaces and codimension-$1$ foliations} \label{ss:foliation}

Here we describe the correspondence between norm-minimizing surfaces and certain codimension-$1$ foliations.

A \textit{codimension-$1$ foliation} $\mathcal{F}$ on a $3$-manifold $M$ is a decomposition of $M$ into possibly noncompact immersed $2$-dimensional submanifolds called leaves such that $M$ is covered by a collection of charts of the form $\R^2 \times \R$ where the leaves pass through a given chart in slices of the form $\R^2 \times \{ z \}$ for $z \in \R$.
A codimension-$1$ foliation $\mathcal{F}$ on $M$ is \textit{transversely oriented} if some vector field on $M$ transverse to the leaves of $\mathcal{F}$ is fixed. 

A Reeb component is the foliation on the solid torus $D^2 \times S^1$ described as follows:
Consider the  decomposition $\R^2 \times [0, \infty)$ into planes $\R^2 \times \{ z \}$ for $z \in [0, \infty)$ and the action $(x, y, z) \to 2 (x, y, z)$ on $\R^2 \times [0, \infty)$ by $\Z$.
The quotient $(\R^2 \times [0, \infty) \setminus \{ (0, 0, 0) \}) / \Z$ is a solid torus, and the induced foliation on it is a Reeb foliation.

As a generalization of the fact that every fiber of a fibration over a circle is norm-minimizing, Thurston~\cite[Corollary 2]{Th86} proved the following theorem.
A $3$-manifold is \textit{irreducible} if every embedded $2$-sphere in $M$ bounds an embedded $3$-ball in $M$.
\begin{thm}[\cite{Th86}] \label{thm:foliation1}
Let $M$ be an irreducible $3$-manifold with empty or toroidal boundary and $\mathcal{F}$ a codimension-$1$ transversely oriented foliation on $M$ such that $\mathcal{F}$ has no Reeb components, and each component of $\partial M$ is either transverse to $\mathcal{F}$ or is a leaf of $\mathcal{F}$.
Then every compact leaf of $\mathcal{F}$ is norm-minimizing.
\end{thm}

For such a codimension-$1$ foliation $\mathcal{F}$ as in Theorem~\ref{thm:foliation1} we can consider the Euler class of the bundle of planes tangent to the leaves.
The key ingredient of the proof is that the dual Thurston norm of the Poincar\'e dual of the Euler class is less than or equal to $1$.

A codimension-$1$ transversely oriented foliation $\mathcal{F}$ on $M$ is \textit{taut} if there exists a closed curve in $M$ transversally intersecting each leaf of $\mathcal{F}$. 
Every taut foliation has no Reeb components~\cite{Go75}.
Let $\mathcal{F}$ be a codimension-$1$ foliation on $M$.
We say that a leaf $L$ of $\mathcal{F}$ is of depth $0$ if $L$ is compact.
Having defined depth $j \leq k$ leaves we say that $L$ is of depth $k + 1$ if $\overline{L} \setminus L$ is a union of depth $j \leq k$ leaves and contains a leaf of depth $k$.   
We say that $\mathcal{F}$ is \textit{of finite depth} if there exists an integer $k$ such that the depth of every leaf of $\mathcal{F}$ is defined to be less than $k$.

Developing the theory of sutured manifolds, Gabai~\cite[Theorem 5.5]{Ga83} proved the following theorem, which can be seen as the converse of Theorem~\ref{thm:foliation1}:
\begin{thm}[\cite{Ga83}] \label{thm:foliation2}
Let $M$ be an irreducible $3$-manifold with empty or toroidal boundary and $S$ a norm-minimizing surface in $M$ representing a nontrivial class in $H_2(M, \partial M; \Z)$.
Then there exists a codimension-$1$ transversely oriented taut foliation $\mathcal{F}$ on $M$ of finite depth such that $\mathcal{F}$ is transverse to $\partial M$, $S$ is a leaf of $\mathcal{F}$ and $\mathcal{F}|_{\partial M}$ is a suspension of homeomorphisms of $S^1$.
\end{thm}

Gabai's construction of such a foliation as in Theorem~\ref{thm:foliation2} used a so-called sutured manifold hierarchy~\cite{Ga83, Sc89}:
A \textit{sutured manifold} $(M, \gamma)$ is a $3$-manifold $M$ equipped with a decomposition of $\partial M$ into two subsurfaces $R_\pm$ meeting along a possibly empty system $\gamma$ of simple closed curves such that $R_-$ and $R_-$ are transversely oriented inwards and outwards respectively.
Under the assumptions in Theorem~\ref{thm:foliation2} there exists a sequence of sutured manifolds $(M_0, \gamma_0), \dots, (M_n, \gamma_n)$, where $(M_0, \gamma_0)$ is obtained by decomposing $M$ along $S$, $(M_i, \gamma_i)$ is obtained by decomposing $M_{i-1}$ along certain type of properly embedded surface, and $(M_n, \gamma_n)$ is a collection of $3$-balls with single simple closed curves.
A codimension-$1$ transversely oriented taut foliation on $M$ is constructed inductively on the hierarchy.

Together with Theorem~\ref{thm:foliation1}, the construction provides an effective algorithm to compute the Thurston norm.
Applying the techniques to knots and links, Gabai~\cite{Ga84} gave tables of the genera of knots with $10$ or fewer crossings and links with $9$ or fewer crossings.

\begin{rem}
Based on the idea, Lackenby~\cite{La21} showed that it is in NP to determine the Thurston norm of a given first cohomology class.
Computations of the Thurston norm were described by Oertel~\cite{Oe86} in terms of branched surfaces, and by Mosher~\cite{Mos91, Mos92a, Mos92b} in terms of pseudo-Anosov flows.
We will see another algorithm via normal surface theory in Section~\ref{ss:normal_surface}. 
\end{rem}

As a corollary of Theorems~\ref{thm:foliation1}, \ref{thm:foliation2}, Gabai~\cite[Corollary 6.13]{Ga83} proved the following, which has been conjectured by Thurston~\cite{Th86}.
\begin{cor}
Let $p \colon \widetilde{M} \to M$ be an $n$-fold covering.
Then
\[ x_{\widetilde{M}}(p^*(\phi)) = n x_M(\phi) \]
for $\phi \in H^1(M; \R)$.
\end{cor}

It is worth pointing out here that for a finite covering $p \colon \widetilde{M} \to M$, $\phi \in H^1(M; \Z)$ is fibered if and only if $p^*(\phi)$ is fibered, which is an immediate consequence of Stallings' fibration theorem.

In \cite[Theorem 2]{Ga87} Gabai showed the following stronger result than Theorem~\ref{thm:foliation2} in the case of knots in $S^3$:

\begin{thm}[\cite{Ga87}] \label{thm:knot-foliation}
Let $K$ be a knot in $S^3$ and $S$ a minimal genus Seifert surface of $K$.
Then there exits a codimension-$1$ taut foliation $\mathcal{F}$ on $K(0)$ of finite depth such that the capped off surface $S$ is a leaf of $\mathcal{F}$.
\end{thm}

As a corollary Gabai~\cite[Corollary 5]{Ga87} proved the following:

\begin{cor}[\cite{Ga87}] \label{cor:zero-surgery}
For a knot $K$ in $S^3$, $K(0)$ is prime and the genus $g(K)$ is equal to the minimal genus of an embedded nonseparating surface in $K(0)$.
\end{cor}

The \textit{Property $R$ conjecture} asserts that if $K(0)$ is homeomorphic to $S^2 \times S^1$, then $K$ is the unknot.
The \textit{Po\'enaru conjecture} is stronger and asserts that if $K(0)$ is reducible, then $K$ is the unknot.
Corollary~\ref{cor:zero-surgery} gave the positive proofs of these conjectures.

\begin{rem}
Let $M$ be an irreducible $3$-manifold with toroidal boundary, not being a cable space and not homeomorphic to $T^2 \times [0, 1]$.
Generalizing a result of Sela~\cite{Sel90}, Baker and Taylor~\cite{BT19} showed that for all but finitely many slopes of $\partial M$, the Thurston norm of $M$ equals that of the result of the Dehn filling along a slope plus the so-called winding norm.
\end{rem}

Another corollary~\cite[Corollary 6]{Ga87} of Theorem~\ref{thm:knot-foliation} is the following:
\begin{cor}[\cite{Ga87}]
A knot $K$ in $S^3$ is fibered if and only if $K(0)$ is fibered.
\end{cor}

Scharlemann~\cite{Sc89} realized that much of Gabai's theory could work only in terms of sutured manifolds without any reference to foliations.
See also \cite{CantC13} for results on a generalization of the Thurston norm for sutured manifolds.

\subsection{Singular and Gromov norms}

Using Theorem~\ref{thm:foliation2}, Gabai~\cite[Corollary 6.18]{Ga83} showed the equivalence of the Thurston norm, its singular one and the Gromov norm of a $3$-manifold $M$.

The \textit{singular Thurston norm} $x_{M, s}$ on $H^1(M; \Z)$ is defined by
\[
\begin{split}
x_{M, s}(\phi) = \min \left\{ \frac{1}{n} \chi_-(S) ~;~ \text{$f \colon (S, \partial S) \to (M, \partial M)$ is a proper map from} \right. \\
\left. \text{a surface $S$ such that $f_*([S, \partial S])$ is dual to $n \phi$} \right\},
\end{split}
\]
for $\phi \in H^1(M; \Z)$.
It is straightforward to see that $x_{M, s}$ uniquely extends to $H^1(M; \R)$ as a seminorm.

For a singular $k$-chain $\sum_i a_i \sigma_i \in C_k(M, \partial M; \R)$ its norm is defined to be the sum $\sum_i |a_i|$ of its absolute values of the coefficients.
The \textit{Gromov norm} or \textit{$l^1$-seminorm} $||c||_1$ of $c \in H_k(M, \partial M; \R)$ is the induced seminorm~\cite{Gro82}:
\[ ||c||_1 = \inf \left\{ \sum_i |a_i| ~;~ \text{$\sum_i a_i \sigma_i$ is a singular $2$-cycle representing $c$} \right\}. \]
We also denote by $||\cdot||_1$ the seminorm on $H^*(M; \R)$ induced by Poincar\'e duality.

The first equality in the following theorem has been conjectured by Thurston~\cite{Th86}.
See also \cite{Pe93} for a combinatorial proof.
\begin{thm}[\cite{Ga83}] \label{thm:Gromov}
The following equality holds on $H^1(M; \R)$:
\[ x_M = x_{M, s} = \frac{1}{2} ||\cdot||_1. \]
\end{thm}

As a special case, Gabai~\cite[Corollary 6.23]{Ga83} proved the following generalization of Dehn's lemma for higher genus surfaces.
\begin{cor}[\cite{Ga83}] \label{cor:Dehn}
Let $f \colon S \to M$ be a map from a surface with connected boundary such that $f|_{\partial S}$ is an embedding and $f^{-1}(f(\partial S)) = \partial S$.
Then there exists an embedded surface $S'$ in $M$ such that $\partial S' = \partial S$ and the genus of $S'$ is less than or equal to that of $S$.
\end{cor}

In particular, Corollary~\ref{cor:Dehn} shows equality of the embedded and immersed genera of knots~\cite[Corollary 6.22]{Ga83}: 
\begin{cor}[\cite{Ga83}]
The genus $g(K)$ of a knot $K$ in $S^3$ is equal to the minimal genus of immersed surfaces $S$ in $S^3$ bounding $K$ which are nonsingular along $K$.
\end{cor}

%%%%%%%% Section 3 %%%%%%%%%%%%%%%%%%%%%%%%%%%%%%%%%%%%%%%%%%%%%%%%%%%%%%%%%%%%
\section{Alexander and Teichm\"uller polynomials} \label{sec:3}

We describe the lower bound on the Thurston norm by the Alexander polynomial, following McMullen~\cite{Mc02}.
This lower bound is then restated in terms of abelian Reidemeister torsion.
We also discuss the Teichm\"uller polynomial associated with a fibered face of the Thurston norm ball, introduced by McMullen~\cite{Mc00}.

\subsection{Alexander polynomial} \label{ss:Alexander_polynomial}

It is well known that the classical Alexander polynomial $\Delta_K(t) \in \Z[t, t^{-1}]$ of a knot $K$ in $S^3$ gives a lower bound on the genus $g(K)$:
 \[ 2 g(K) \geq \deg \Delta_K(t), \]
where equality holds if $K$ is a fibered knot.
Following McMullen~\cite{Mc02}, we describe a generalization of this inequality on the Thurston norm and the Alexander polynomial of a general $3$-manifold.

Let $M$ be a $3$-manifold with empty or toroidal boundary.
We denote by $H_1(M)_f$ the free abelian group obtained by dividing $H_1(M; \Z)$ by the torsion submodule.
We denote by $\overline{M}$ the maximal free abelian cover of $M$, which is the cover of $M$ associated with the canonical projection $\pi_1 M \to H_1(M)_f$.
Since $H_1(M)_f$ acts on $\overline{M}$ by deck transformations, $H_1(\overline{M}; \Z)$ is a finitely generated module over the group ring $\Z[H_1(M)_f]$.
The \textit{Alexander polynomial} $\Delta_M \in \Z[H_1(M)_f]$ of $M$ is the order of  $H_1(\overline{M}; \Z)$ over $\Z[H_1(M)_f]$, which is well-defined up to multiplication by elements of $\pm H_1(M)_f$:
In general, for a finitely generated module $L$ over a noetherian UFD $R$ and an exact sequence
\[ R^l \xrightarrow{r} R^m \to L \to 0 \] 
with $l \geq m$, the \textit{order} of $L$ is the greatest common divisor of the $m$-minors of a representation matrix $r$, and is well-defined up to multiplication by units in $R$.

\begin{exmp}
For a knot $K$ in $S^3$, $\Delta_{X_K}$ coincides with the classical Alexander polynomial $\Delta_K(t) \in \Z[t, t^{-1}]$ under the identification of $H_1(X_K; \Z)$ with the infinite cyclic group generated by $t$.
\end{exmp}

McMullen~\cite[Theorem 1.1]{Mc02} introduced the \textit{Alexander norm} $||\cdot||_A$ on $H^1(M; \R)$ and showed an inequality between the Thurston and Alexander norms as follows:
We write $\Delta_M = \sum_{h \in H_1(M)_f} a_h h$ for $a_h \in \Z$.
If $\Delta_M = 0$, then we define $||\cdot||_A = 0$.
Otherwise, we define
\[ ||\phi||_A = \max \{ \langle \phi, h - h' \rangle ~;~ h, h' \in H_1(M)_f ~\text{such that}~ a_h a_{h'} \neq 0 \} \]
for $\phi \in H^1(M; \R)$.
It is clear that $||\cdot||_A$ is a seminorm on $H^1(M; \R)$.

\begin{thm}[\cite{Mc02}] \label{thm:A-norm}
Let $M$ be a $3$-manifold with empty or toroidal boundary.
Then
\[ x_M(\phi) \geq ||\phi||_A -
\begin{cases}
1 + b_3(M) &\text{if $H^1(M; \Z)$ is generated by $\phi$}, \\
0 &\text{if $b_1(M) > 1$},
\end{cases}
\]
for $\phi \in H^1(M; \Z)$.
Furthermore, equality holds if $\phi$ is fibered with $M \neq S^1 \times S^2$ and $M \neq S^1 \times D^2$.
\end{thm}

\begin{rem}
As Dunfield~\cite{Du01} showed, there are examples of $3$-manifolds $M$ fibering over a circle with $b_1(M) > 1$ such that $x_M$ and $||\cdot||_A$ do not agree.
\end{rem}

It is known that $2 g(K) = \deg \Delta_K(t)$ for all knots up to $10$ crossings or less (See for example \cite{Ga84}). 
McMullen~\cite[Theorem 7.1]{Mc02} showed that the Thurston and Alexander norms agree for all the tabulated links with $9$ or fewer crossings in \cite{Ro76} except $9_{21}^3$, and possibly $9_{41}^2$, $9_{50}^2$ and $9_{15}^3$.

\subsection{Abelian torsion}

We discuss a corresponding result to Theorem~\ref{thm:A-norm} in terms of abelian Reidemeister torsion.
We will see the precise definition of Reidemeister torsion in Section~\ref{ss:R-torsion}.

Let $M$ be a $3$-manifold with empty or toroidal boundary with a CW-complex structure.
We denote by $\Q(H_1(M)_f)$ the quotient field of $\Z[H_1(M)_f]$.
The \textit{abelian torsion} or \textit{Milnor torsion} $\tau(M) \in \Q(H_1(M)_f)$ of $M$ is the Reidemeister torsion associated with the canonical projection $\pi_1 M \to H_1(M)_f$, which is the algebraic torsion of the twisted chain complex $C_*(\overline{M}) \otimes_{\Z[H_1(M)_f]} \Q(H_1(M)_f)$ of the CW-complex $\overline{M}$.
The topological invariant $\tau(M)$ is well-defined up to multiplication by elements of $\pm H_1(M)_f$, and is known to be symmetric, i.e., $\tau(M)$ is invariant up to multiplication by elements of $\pm H_1(M)_f$ under the involution on $\Q(H_1(M)_f)$ reversing the elements of $H_1(M)_f$.

\begin{exmp}
For a knot $K$ in $S^3$, 
\[ \tau(X_K) = \frac{\Delta_K(t)}{t-1} \]
under the identification of $H_1(X_K; \Z)$ with the infinite cyclic group generated by $t$.
\end{exmp}

Turaev~\cite{Tu75, Tu02b} showed that $\tau(M)$ determines the Alexander polynomial $\Delta_M$, and vise versa:
\begin{thm}[\cite{Tu75, Tu02b}] \label{thm:abelian_torsion}
Let $M$ be a $3$-manifold with empty or toroidal boundary.
If $H_1(M)_f$ is an infinite cyclic group generated by $t$, then
\[ \tau(M) =
\begin{cases}
\frac{\Delta_M}{(t-1)^2} &\text{if $\partial M = \emptyset$}, \\
\frac{\Delta_M}{t-1} &\text{if $\partial M \neq \emptyset$}.
\end{cases}
\]
If $b_1(M) > 1$, then
\[ \tau(M) = \Delta_M. \]
\end{thm}

A cohomology class $\phi \in H^1(M; \Z)$ induces a ring homomorphism $\Z[H_1(M)_f] \to \Z[t, t^{-1}]$ by sending $h \in H_1(M)_f$ to $t^{\langle \phi, h \rangle}$.
We define $\tau_\phi(M) \in \Q(t)$ to be the reduction of $\tau(M)$ by the induced homomorphism, which is the algebraic torsion of $C_*(\overline{M}) \otimes_{\Z[H_1(M)_f]} \Q(t)$.

Theorem~\ref{thm:A-norm} is restated in terms of $\tau_\phi$ as follows.
We define
\[\deg \left( a_l t^l + a_{l+1} t^{l+1} + \cdots + a_m t^m \right) = m - l \]
for $l$, $m \in \Z$ with $l < m$ and $a_i \in \Z$ with $a_l a_m \neq 0$, and further define
\[ \deg \frac{p(t)}{q(t)} = \deg p(t) - \deg q(t) \]
for $p(t), q(t) \in \Z[t, t^{-1}] \setminus \{ 0 \}$.

\begin{thm} \label{thm:torsion-norm}
Let $M$ be a $3$-manifold with empty or toroidal boundary.
Then
\[ x_M(\phi) \geq \deg \tau_\phi(M) \]
for $\phi \in H^1(M; \Z)$.
Furthermore, equality holds if $\phi$ is fibered with $M \neq S^1 \times S^2$ and $M \neq S^1 \times D^2$.
\end{thm}

It is well known that $\Delta_K(t)$ is monic for a fibered knot $K$ in $S^3$.
More generally, $\tau_\phi(M)$ is represented by a monic polynomial divided by $(t-1)^{1 + b_3(M)}$ for a fibered class $\phi \in H^1(M; \Z)$.
We will discuss more the property in Remark~\ref{rem:monic}.

In analogy with the Thurston norm, Turaev~\cite{Tu02a, Tu07} introduced a seminorm on $H^1(X; \R)$ for a finite $2$-dimensional complex $X$, and numerical functions on $H_2(M; \Q / \Z)$ and on the torsion subgroup of $H_1(M; \Z)$ for a $3$-manifold $M$.
Turaev showed that the Alexander polynomial and abelian Reidemeister torsion give lower bounds also on these functions. 
See \cite{FSW16, NW14} for further studies of such analogues of the Thurston norm.

\subsection{Teichm\"uller polynomial}

Here we give a brief exposition of the Teichm\"uller polynomial introduced by McMullen~\cite{Mc00}.

A (codimension-$1$) lamination $\mathcal{L}$ on a $3$-manifold $M$ is a codimension-$1$ foliation on a closed subset of $M$.
A lamination $\mathcal{L}$ on $M$ is \textit{transversely orientable} if there is a nonsingular vector field on a neighborhood of the underlying closed subset of $\mathcal{L}$ in $M$ transverse to the leaves.
A \textit{geodesic lamination} $\lambda$ on a hyperbolic surface $S$ is a decomposition of a closed subset of $S$ into simple geodesics.

Let $M$ be a hyperbolic $3$-manifold having a fibered face $F$ of $B_M$.
Let $\phi \in H^1(M; \Z)$ be a fibered class in the cone on $F$, and $\psi \colon S \to S$ a pseudo-Anosov monodromy of a fibration representing $\phi$.
Then $\psi$ has an expanding invariant geodesic lamination $\lambda$ on $S$.
Let $\mathcal{L}$ be the lamination on $M$ obtained as the mapping torus of $\psi|_\lambda$.
Based on results by Fried~\cite{Frid82a}, McMullen~\cite[Corollary 3.2]{Mc00} showed that the isotopy class of $\mathcal{L}$ depends only on $F$.

We denote by $\overline{\mathcal{L}}$ the preimage of $\mathcal{L}$ by the maximal free abelian covering $\overline{M} \to M$.
A \textit{transeversal} for $\overline{\mathcal{L}}$ is a compact totally disconnected subset of $\overline{\mathcal{L}}$ such that there is an open neighborhood $U$ of $T$ with a homeomorphism $(U, T) \to (T \times \R^2, T \times \{ 0 \})$.
Note that the free abelian group $H_1(M)_f$ acts on the set of transversals for $\mathcal{L}$.
We define $T(\overline{\mathcal{L}})$ to be the abelian group generated by all transversals $[T]$ for $\overline{\mathcal{L}}$ modulo the following relations:
\begin{enumerate}
\item $[T] = [T'] + [T'']$, ~if $T$ is a disjoint union of $T'$ and $T''$,
\item $[T] = [T']$, ~if there is an open neighborhood $U$ of $T \cup T'$ with homeomorphisms $(U, T) \to (T \times \R^2, T \times \{ 0 \})$, $(U, T') \to (T' \times \R^2, T' \times \{ 0 \})$.
\end{enumerate}
A consequence of the compactness of $\mathcal{L}$ is that $T(\overline{\mathcal{L}})$ is a finitely generated $\Z[H_1(M)_f]$-module.

Now the \textit{Teichm\"uller polynomial} $\Theta_F \in \Z[H_1(M)_f]$ of $F$ is defined to be the order of $T(\overline{\mathcal{L}})$ over $\Z[H_1(M)_f]$, which is well-defined up to multiplication by elements of $\pm H_1(M)_f$.
McMullen showed that $\Theta_F$ is monic and symmetric.

McMullen~\cite[Theorem 6.1]{Mc00} introduced the \textit{Teichm\"uller norm} $||\cdot||_{\Theta_F}$ on $H^1(M; \R)$ and showed its relation with the Thurston norm as follows:
We write $\Theta_F = \sum_{h \in H_1(M)_f} a_h h$ for $a_h \in \Z$ and define
\[ ||\phi||_{\Theta_F} = \max \{ \langle \phi, h - h' \rangle ~;~ h, h' \in H_1(M)_f ~\text{such that}~ a_h a_{h'} \neq 0 \} \]  
for $\phi \in H^1(M; \R)$.
It is clear that $||\cdot||_{\Theta_F}$ is a seminorm on $H^1(M; \R)$.

\begin{thm}[\cite{Mc00}] \label{thm:T-norm}
Let $F$ be a fibered face of $B_M$ of a hyperbolic $3$-manifold $M$.
Then there exists a face $D$ of the unit ball of $||\cdot||_{\Theta_F}$ such that the cones on $F$ and $D$ coincides.  
\end{thm}

Together with a computational formula of $\Theta_F$ in terms of train tracks on fibers~\cite[Theorem 3.6]{Mc00}, Theorem~\ref{thm:T-norm} provides an effective algorithm to determine a fibered face of $B_M$ for a hyperbolic $3$-manifold $M$ from a single fiber and the monodromy on it.

Using Theorem~\ref{thm:A-norm}, McMullen~\cite[Theorem 7.1]{Mc00} proved the following theorem:
\begin{thm}[\cite{Mc00}]
Let $F$ be a fibered face of $B_M$ of a hyperbolic $3$-manifold $M$ with $b_1(M) > 1$.
Then there exists a unique face $A$ of the unit ball of the Alexander norm containing $F$.
Furthermore, if the lamination $\mathcal{L}$ associated with $F$ is transversely orientable, then $F = A$ and $\Delta_M$ divides $\Theta_F$.
\end{thm}

McMullen~\cite{Mc00} also showed that for a fibered class $\phi \in H^1(M; \Z)$ lying in the cone on a fibered face $F$, the dilatation $\lambda(\phi)$ of its monodromy is the largest root of the polynomial equation $\sum_{h \in H} a_h t^{\phi(h)} = 0$ obtained by evaluating $\Theta_F$ by $\phi$, and that the function $\frac{1}{\log \lambda(\phi)}$ extends to the cone on $F$ as a real-analytic function which is strictly concave, extending results in \cite{Frid82b, Mats87}.
See also \cite{Su15}.

Dowdall, Kapovich and Leininger~\cite{DKL15, DKL17} introduced analogues of the Teichm\"uller polynomial and proved analogous results for free-by-cyclic groups.
In their work a hyperbolic $3$-manifold fibering over a circle and its fibered face of the Thurston norm ball are replaced by a free-by-cyclic group and a component of its Bieri--Neumann--Strebel invariant~\cite{BNS87}.

%%%%%%%% Section 4 %%%%%%%%%%%%%%%%%%%%%%%%%%%%%%%%%%%%%%%%%%%%%%%%%%%%%%%%%%%%
\section{Seiberg--Witten invariant} \label{sec:4}

Here we are concerned with adjunction inequalities, which give relationships between the Seiberg--Witten invariant of a $4$-manifold and the complexity of embedded surfaces in the manifold, and between the Seiberg--Witten invariant of a $3$-manifold and its Thurston norm.
As a related topic we also discuss the harmonic norm on the cohomology group associated with a Riemannian metric.

\subsection{Seiberg--Witten theory} \label{ss:SW}

We briefly review Seiberg--Witten theory~\cite{Wit94} in the case of a closed smooth $4$-manifold with $b_2^+(N) > 1$.
(Here $b_2^+(N)$ is the dimension of a maximal positive-definite subspace $H_+^2(N; \R)$ of the intersection pairing on $H^2(N; \R)$.)
For the details we refer the reader to the expositions \cite{HT99, KrMr07, Linf16, Morg96}.

Recall that the Lie group $\Spin^c(n) = \Spin(n) \times_{\pm 1} U(1)$ is a central extension of $SO(n)$ by $U(1)$.
A $\Spin^c$-structure on a Riemannian $n$-manifold $X$ is a lifting of the principal $SO(n)$-frame bundle on $X$ to a principal $\Spin^c(n)$-bundle.
We deonte by $\Spin^c(X)$ the set of equivalence classes of $\Spin^c$-structures on $X$.
The set $\Spin^c(X)$ has a free and transitive action by $H^2(X; \Z)$, and we write $\mathfrak{s} + c$ for the image of $\mathfrak{s}$ by $c \in H^2(X; \Z)$.
The first Chern class of the principal $U(1)$-bundle associated with a $\Spin^c$-structure on $X$ defines the map $c_1 \colon \Spin^c(X) \to H^2(X; \Z)$.

Let $N$ be a closed Riemannian $4$-manifold with a metric $g$, and $\widetilde{P}$ a spin$^c$ structure on $N$.
We denote by $L$ the determinant line bundle of $\widetilde{P}$ and by $S^\pm$ the two complex spin bundles associated with $\widetilde{P}$. 
For a connection on $L$ we have Dirac operators $D_A \colon \Gamma(S^\pm) \to \Gamma(S^\mp)$ on the set of sections of $S^\pm$, defined using Levi--Civit\'a connection on the frame bundle on $X$.
The \textit{Seiberg--Witten monopole equations} associated with $\widetilde{P}$ are the following pair of nonlinear elliptic equations for unitary connections $A$ on $L$ and sections $\psi$ of $S^+$:
\begin{align*}
F_A^+ &= \psi \otimes \psi^* - \frac{|\psi|^2}{2} Id, \\
D_A(\psi) &= 0,
\end{align*}
where we identify $S^+$ with its dual via an anti-complex isomorphism and $\psi^*$ is its image of $\psi$.
We denote by $\mathcal{M}_N(\widetilde{P})$ the quotient of the space of gauge-equivalence classes of solutions to the equations.
The moduli space $\mathcal{M}_N(\widetilde{P})$ is known to be compact~\cite{Wit94}.
A class $c \in H^2(N; \Z)$ is called a \textit{Seiberg--Witten monopole class} if $c = c_1(\mathfrak{s})$ for some $\mathfrak{s} \in \Spin^c(M)$ representing $\widetilde{P}$ with nonempty $\mathcal{M}(\widetilde{P}, g)$.

If $b_2^+(N) > 1$, then the \textit{Seiberg--Witten invariant} $SW_N \colon \Spin^c(N) \to \Z$ is defined as follows:
Let $\widetilde{P}$ be a $\Spin^c$-structure on $N$ representing $\mathfrak{s} \in \Spin^c(N)$.
We denote by $\mathcal{C}_N(\widetilde{P})$ the space of gauge-equivalence classes of pairs $(A, \psi)$ with $\psi \neq 0$, which is a classifying space of the group $(S^1)^N$.
There is a universal $S^1$-bundle over $\mathcal{C}_N(\widetilde{P})$ whose Chern class $\mu$ generates $H^2(\mathcal{C}_N(\widetilde{P}); \Z)$.
After a perturbation of $\mathcal{M}_N(\widetilde{P})$ by an addition of $i \eta$ to the right hand side of the first Seiberg--Witten monopole equation for a generic (real) self-dual $2$-form $\eta$ on $N$, the resulting moduli space $\mathcal{M}$ is known to become a smooth submanifold in $\mathcal{C}_N(\widetilde{P})$ of dimension
\[ d(\mathfrak{s}) = \frac{\langle c_1(\mathfrak{s})^2, [N] \rangle - (2 \chi(N) + 3 \sigma(N))}{4}, \]
where $\sigma(N)$ is the signature of $N$.
Moreover, choosing an orientation of the real vector space $H_+^2(N; \R) \oplus H^1(N; \R)$ gives an orientation of $\mathcal{M}$.
If $d(\mathfrak{s})$ is odd, then we define $SW_N(\mathfrak{s}) = 0$, and otherwise we define
\[ SW_N(\mathfrak{s}) = \langle \mu^{d(\mathfrak{s})}, [\mathcal{M}] \rangle. \]
This is an invariant of $\mathfrak{s}$, which is independent of the choice of the Riemannian metric of $N$ and the perturbation term $\eta$~\cite{Wit94}.
The invariant $SW_N$ takes nonzero value only on finitely many $\Spin^c$-structures, and changing the orientation of $H_+^2(N; \R) \oplus H^1(N; \R)$ reverses its sign.
A class $c \in H^2(N; \Z)$ is called a \textit{Seiberg--Witten basic class} if $c = c_1(\mathfrak{s})$ for some $\mathfrak{s} \in \Spin^c(N)$ with $SW_N(\mathfrak{s}) \neq 0$.
Note that every Seiberg--Witten basic class is a Seiberg--Witten monopole class.

\subsection{Seiberg--Witten invariant of a $3$-manifold}

The \textit{Seiberg--Witten invariant} $SW_M \colon \Spin^c(M) \to \Z$ of a closed $3$-manifold $M$ with $b_1(M) > 1$ can be defined by 
\[ SW_M(\mathfrak{s}) = SW_{M \times S^1}(\pi^* \mathfrak{s}) \]
for $\mathfrak{s} \in \Spin^c(M)$, where $\pi^* \mathfrak{s} \in \Spin^c(M \times S^1)$ is the pullback of $\mathfrak{s}$ by the projection $\pi \colon M \times S^1 \to M$.
Note that all solutions to the Seiberg--Witten monopole equations on $M \times S^1$ are known to be $S^1$-invariant~\cite{OT98}.
One can also define $SW_M$ directly in terms of Seiberg--Witten monopole equations on $M$ as in Section~\ref{ss:SW}.
As in the case of a $4$-manifold, a class $c \in H^2(M; \Z)$ is called a \textit{Seiberg--Witten monopole class} if $c = c_1(\mathfrak{s})$ for some $\mathfrak{s} \in \Spin^c(M)$ representing a $\Spin^c$-structure $\widetilde{P}$ on $M$ with nonempty $\mathcal{M}_{M \times S^1}(\pi^* \widetilde{P})$ for the pullback $\pi^* \widetilde{P}$ on $M \times S^1$.
Also, a class $c \in H^2(M; \Z)$ is called a \textit{Seiberg--Witten basic class} if $c = c_1(\mathfrak{s})$ for some $\mathfrak{s} \in \Spin^c(M)$ with $SW_M(\mathfrak{s}) \neq 0$.

For a closed $3$-manifold $M$, Turaev~\cite{Tu89, Tu97} introduced a refinement of the abelian torsion $\tau(M)$ as an integer-valued function $T_M$ on $\Spin^c(M)$ (or on the set of so-called Euler structures on $M$), called \textit{Turaev's torsion function} of $M$.
When $b_1(M) > 1$, $\tau(M)$ is represented by an element of $\Z[H_1(M)_f]$ (Theorem~\ref{thm:abelian_torsion}), and $T_M$ satisfies
\[ \tau(M) = \sum_{h \in H_1(M; \Z)} T_M(\mathfrak{s} - PD(h)) [h] \]
for $\mathfrak{s} \in \Spin^c(M)$, where $PD(h) \in H^2(M; \Z)$ is the Poincar\'e dual of $h$.
Similarly, when $H_1(M)_f$ is an infinite cyclic group generated by $t$, $\tau(M)$ is represented by an element of the Novikov ring $\Z((t)) = \Z[[t]][t^{-1}]$ (Theorem~\ref{thm:abelian_torsion}), and $T_M$ satisfies
\[ \tau(M) = \sum_{h \in H_1(M; \Z)} T_M(\mathfrak{s} - PD(h)) [h] \in \Z((t)) \]
for $\mathfrak{s} \in \Spin^c(M)$.

In terms of $T_M$, Turaev~\cite[Theorem 1]{Tu98} refined the equivalence of the Seiberg--Witten invariant and the abelian torsion shown by Meng and Taubes~\cite{MeTa96}:
\begin{thm}[\cite{MeTa96, Tu98}] \label{thm:MT}
For a closed $3$-manifold $M$ with $b_1(M) > 1$, the Seiberg--Witten invariant and Turaev's torsion function of $M$ coincides up to sign:
\[ SW_M = \pm T_M. \]
\end{thm}

Theorem~\ref{thm:MT} similarly extends to the case $b_1(M) = 1$~\cite{MeTa96, Tu98}.
See \cite{Nic04} for the case $b_1(M) = 0$.

The following is the \textit{adjunction inequality} for $3$-manifolds.
See \cite{Au96, Kr98} for the details.
We will discuss more on adjunction inequalities in Section~\ref{ss:adjunction}.

\begin{thm}[\cite{Au96, Kr98}] \label{thm:basic}
Let $M$ be a closed irreducible $3$-manifold with $b_1(M) > 1$ and $c \in H^2(M; \Z)$ a Seiberg--Witten basic class.
Then
\[ x_M(\phi) \geq |\langle c \cup \phi, [M] \rangle| \]
for $\phi \in H^1(M; \R)$.  
\end{thm}

Kronheimer and Mrowka~\cite[Theorem 1]{KrMr97b} showed that the Thurston norm is determined by the Seiberg--Witten monopole classes:

\begin{thm}[\cite{KrMr97b}] \label{thm:monopole}
Let $M$ be a closed irreducible $3$-manifold with $b_1(M) > 1$.
Then 
\[ x_M(\phi) = \max \{ |\langle c \cup \phi, [M] \rangle| ~;~ \text{$c \in H^2(M; \Z)$ is a Seiberg--Witten monopole class} \} \]
for $\phi \in H^1(M; \R)$.
\end{thm}

\begin{cor}[\cite{KrMr97b}] 
Let $M$ be a closed irreducible $3$-manifold with $b_1(M) > 1$.
Then the convex hull of the Seiberg--Witten monopole classes in $H^2(M; \R)$ is equal to $B_M^*$. 
\end{cor}

As described by Kronheimer~\cite{Kr98} and Vidussi~\cite{V03}, Theorems~\ref{thm:MT}, \ref{thm:basic}, \ref{thm:monopole} deduce Theorem~\ref{thm:A-norm} for closed irreducible $3$-manifolds $M$ with $b_1(M) > 1$.

\subsection{Complexity of surfaces in a $4$-manifold} \label{ss:adjunction}

Adjunction inequalities give relationships between the Seiberg--Witten invariants of a $4$-manifold and the genus of embedded surfaces in the manifold.
The terminology arises from the adjunction formula for a smooth algebraic curve $C$ in an algebraic surface $X$:
\[ \chi_-(C) = C \cdot C - \langle c_1(X), [C] \rangle. \]

The genus of the algebraic curve of degree $d$ in $\C P^2$ is given by $\frac{(d-1)(d-2)}{2}$.
A conjecture attributed to Thom states that the genus of the algebraic curve is minimal among smoothly embedded surfaces in $\C P^2$ representing $d [\C P^1] \in H_2(\C P^2; \Z)$.
With the advance of the Seiberg--Witten monopole equations, Kronheimer and Mrowka~\cite{KrMr94, KrMr95} and Morgan, Szab\'o and Taubes~\cite{MST96} proved the Thom conjecture for holomorphic curves in a general K\"ahler surface with nonnegative intersection.
Later, Ozsv\'ath and Szab\'o~\cite[Theorem 1.1 and Corollary 1.2]{OS00a} proved the symplectic Thom conjecture in its complete generality:

\begin{thm}[\cite{OS00a}] \label{thm:symplectic}
The genus of an embedded symplectic surface in a closed symplectic $4$-manifold is minimal among smoothly embedded surfaces representing the same homology class. 
\end{thm}

The theorem for K\"ahler surfaces follows as a special case of Theorem~\ref{thm:symplectic}:
\begin{cor}
The genus of an embedded holomorphic curve in a K\"ahler surface is minimal among smoothly embedded surfaces representing the same homology class.
\end{cor}

The following are adjunction inequalities shown by Morgan, Szab\'o and Taubes~\cite[Proposition 4.2]{MST96} and Ozsv\'ath and Szab\'o~\cite[Corollary 1.7]{OS00a}:

\begin{thm}[\cite{MST96}] \label{thm:adjunction1}
Let $N$ be a smooth closed $4$-manifold with $b_2^+(M) > 1$ and $c \in H^2(N; \Z)$ a Seiberg--Witten basic class.
Then for a smoothly embedded surface $\Sigma$ in $N$ with nonpositive Euler charcteristic and $[\Sigma] \cdot [\Sigma] \geq 0$, we have
\[ \chi_-(\Sigma) \geq [\Sigma] \cdot [\Sigma] + \langle c, [\Sigma] \rangle. \]
\end{thm}

\begin{thm}[\cite{OS00a}] \label{thm:adjunction2}
Let $N$ be a smooth closed $4$-manifold with $b_2^+(M) > 1$ and $c \in H^2(N; \Z)$ a Seiberg--Witten basic class.
Suppose that $d(\mathfrak{s}) = 0$ for any $\mathfrak{s} \in \Spin^c(N)$ associated with a basic class in $H^2(N; \Z)$.
Then for a smoothly embedded surface in $N$ with nonpositive Euler charcteristic and $[\Sigma] \cdot [\Sigma] < 0$, we have
\[ \chi_-(\Sigma) \geq [\Sigma] \cdot [\Sigma] + |\langle c, [\Sigma] \rangle|. \] 
\end{thm}

See \cite{OS00b} for further refinements of Theorems~\ref{thm:adjunction1}, \ref{thm:adjunction2}.

There is also an adjunction inequality by Fintushel and Stern~\cite{FS95} for embedded spheres:
\begin{thm}[\cite{FS95}] \label{thm:FS}
Let $N$ be a smooth closed $4$-manifold with $b_2^+(M) > 1$.
Suppose that there exists a Seiberg--Witten basic class.
Then there exist no smoothly embedded spheres $\Sigma$ such that $\Sigma \cdot \Sigma \geq 0$ and $[\Sigma] \neq 0$.
\end{thm}

See \cite{Ba04, BaFu04} for a refinement of Theorem~\ref{thm:FS} in terms of the so-called Bauer-Furuta invariants.
See also \cite{Str03, Ko16} for adjunction-type inequalities for families of embedded surfaces.

Now we describe results on the relationship between complexity of embedded surfaces in circle bundles over a $3$-manifold and its Thurston norm.

Let $N$ be a smooth closed $4$-manifold.
We define a function $x_N \colon H_2(N; \Z) \to \Z$ by
\[ x_N(\alpha) = \min \{ \chi_-(\Sigma) ~;~ \text{$\Sigma$ is an embedded surface representing $\alpha$} \} \]
for $\alpha \in H_2(N; \Z)$.

Using Agol's virtual fibering theorem~\cite{Ag08, Ag13} (see Theorem~\ref{thm:virtual_fibering}) and considering the Seiberg--Witten invariants of finite covers, Friedl and Vidussi~\cite[Theorem 1.1]{FV14b}, and Nagel~\cite[Theorem 5.6]{Nag16} showed the following theorem:
\begin{thm}[\cite{FV14b, Nag16}] \label{thm:complexity1}
Let $M$ be a closed irreducible $3$-manifold which is not a Seifert fibered space and not covered by a torus bundle, and let $p \colon N \to M$ be an oriented circle bundle.
Then 
\[ x_N(\alpha) \geq |\alpha \cdot \alpha| + x_M(p_* \alpha) \]
for $\alpha \in H_2(N; \Z)$.
\end{thm}

\begin{rem}
Kronheimer~\cite{Kr99} proved the same inequality as in Theorem~\ref{thm:complexity1} for the case $N = M \times S^1$ such that $M$ is a closed irreducible $3$-manifold whose Thurston norm does not identically vanish.
\end{rem}

Let $M$ be a closed $3$-manifold.
We set $\Xi_M$ to be the inverse image by the canonical map $H^2(M; \Z) \to H^2(M; \R)$ of the set of nonzero classes $w \in H^2(M; \R)$ such that $v + 2 w$ lies on an edge of $B_M^*$ for some vertex $v$ of $B_M^*$.
Note that $\Xi_M$ is a finite set.

Friedl and Vidussi~\cite[Corollary 1.3]{FV14b} also showed that equality in Theorem~\ref{thm:complexity1} holds for all but finitely many circle bundles over a $3$-manifold which is not exceptional:
\begin{thm}[\cite{FV14b}] \label{thm:complexity2}
Let $M$ be a closed irreducible $3$-manifold which is not a closed graph manifold such that $\Delta_M^\phi \neq 0$ for all nontrivial $\phi \in H^1(M; \Z)$, and let $p \colon N \to M$ be an oriented circle bundle with Euler class not in $\Xi_N$. 
Then
\[ x_N(\alpha) = |\alpha \cdot \alpha| + x_M(p_* \alpha) \]
for $\alpha \in H_2(N; \Z)$.
\end{thm}

\begin{rem}
Friedl and Vidussi showed Theorem~\ref{thm:complexity2} also for nonpositively curved graph manifolds.
By the work by Agol~\cite{Ag13}, Liu~\cite{Liu13}, and Przytycki and Wise~\cite{PW18, Wis21} the so-called virtually special theorem holds for irreducible nonpositively curved $3$-manifolds.
As a consequence, Agol's virtual fibering theorem also holds for such $3$-manifolds.
\end{rem}

\subsection{Harmonic norm}

We discuss relationships between the harmonic norm associated with a Riemannian metric and the Thurston norm. 

Let $M$ be a closed Riemannian $3$-manifold with a metric $h$.
The $L^2$-norm $||\cdot||_h$ on the vector space $\Omega^k(M)$ of $k$-forms on $M$ is associated with the inner product
\[ \langle \alpha, \beta \rangle = \int_M \alpha \wedge * \beta \]
for $\alpha$, $\beta \in \Omega^k(M)$, where $*$ is the Hodge star operator.
The $k$-th homology group $H^k(M; \R)$ is identified with the subspace of harmonic $k$-forms, and $||\cdot||_h$ induces a norm on $H^k(M; \R)$.
The induced norm is called the \textit{harmonic norm} and is also denoted by $||\cdot||_h$.
As it comes from a positive-definite inner product, the unit ball of $||\cdot||_h$ is a smooth ellipsoid.

In the study of the Seiberg--Witten monopole equations Kronheimer and Mrowka~\cite{KrMr97a, KrMr97b} showed that the Thurston norm $x_M$ is characterized  in terms of the harmonic norm:
\begin{thm}[\cite{KrMr97b}] \label{thm:KM}
Let $M$ be  a closed irreducible $3$-manifold not containing non-separating tori.
Then 
\[ x_M(\phi) = \frac{1}{4 \pi} \inf_h ||s_h||_h ||\phi||_h \]
for $\phi \in H^1(M; \R)$, where $s_h$ is the scalar curvature of $h$, and the infimum is taken over all Riemannian metrics $h$ on $M$.
\end{thm}

\begin{rem}
The original statement of \cite[Theorem 2]{KrMr97b} is in terms of the dual Thurston norm $x_M^*$ on $H_2(M; \R)$.
As explained in \cite[Theorem 5.1]{BD17}, it is equivalent to that of Theorem~\ref{thm:KM}.
\end{rem}

See \cite{BS19, Ste19} for extensions of Theorem~\ref{thm:KM} in another approach studying harmonic $1$-forms.
See also \cite{Kat05}.

For a closed hyperbolic $3$-manifold $M$, by Mostow rigidity, the harmonic norm is uniquely determined by the underlying topology of $M$.
We denote it by $||\cdot||_{L^2}$.
Refining results of Bergeron, \c{S}eng\"un and Venkatesh~\cite{BSV16} as well as Theorem~\ref{thm:KM}, Brock and Dunfield~\cite{BD17} showed the following inequalities between the two norms $||\cdot||_{L^2}$ and $x_M$:

\begin{thm}[\cite{BD17}] \label{thm:BD}
Let $M$ be a closed hyperbolic $3$-manifold.
Then
\[ \frac{\pi}{\sqrt{\vol(M)}} x_M(\phi) \leq ||\phi||_{L^2} \leq \frac{10 \pi}{\sqrt{\inj(M)}} x_M(\phi) \]
for $\phi \in H^1(M; \R)$, where $\inj(M)$ is the injectivity radius of $M$, which is half the length of the shortest closed geodesic in $M$.
\end{thm}

Brock and Dunfield used the theory of minimal surfaces to prove Theorem~\ref{thm:BD}.
They also showed that the inequality is qualitatively sharp~\cite[Theorem 1.3]{BD17}.

Since the scalar curvature of a hyperbolic metric $h$ is $-6$, specializing Theorem~\ref{thm:KM} to such $h$ gives
\[ \frac{2 \pi}{3 \sqrt{\vol(M)}} x_M(\phi) \leq ||\phi||_{L^2}, \]
which is weaker than the first inequality in Theorem~\ref{thm:BD}.
Lin~\cite{Linf20} gave a gauge-theoretic proof of the stronger inequality.

%%%%%%%% Section 5 %%%%%%%%%%%%%%%%%%%%%%%%%%%%%%%%%%%%%%%%%%%%%%%%%%%%%%%%%%%%
\section{Floer homology} \label{sec:5}

We look at the fact that Floer homology detects the Thurston norm and fiberedness of a $3$-manifold.

With a motivation to better understand the Seiberg--Witten invariant, Ozsv\'ath and Szab\'o~\cite{OS04d, OS04e} introduced \textit{Heegaard Floer homology}:
\[ \widehat{HF}(M),~ HF^\infty(M),~ HF^+(M),~ HF^-(M). \]
Analogously, based directly on the Seiberg--Witten monopole equations, Kronheimer and Mrowka~\cite{KrMr07} introduced \textit{monopole Floer homology}:
\[ \widecheck{HM}_*(M),~ \widehat{HM}_*(M),~ \overline{HM}_*(M). \]
Passing through \textit{embedded contact homology} $ECH$ introduced by Hutchings and Taubes~\cite{Hu02, HT07, HT09}, Heegaard Floer homology and Monopole Floer homology were shown to be equivalent by Colin, Ghiggini and Honda~\cite{CGH11, CGH12a, CGH12b}, and Kutluhan, Li and Taubes~\cite{KLT20a, KLT20b, KLT20c, KLT20d, KLT20e}.

In the following we focus on results in terms only of Heegaard Floer homology.
But by the equivalence of the theories corresponding results hold also in terms of monopole Floer homology.
For details including the definitions of the Floer homology groups, we refer the reader to the expositions \cite{Gre21, Ho17, Ju15, OS04a, OS04c, OS06, OS18} for Heegaard Floer homology, and to \cite{KrMr07, Linf16} for monopole Floer homology.
For combinatorial computations of Heegaard Floer homology see the survey article~\cite{Man13} and the references given there.

\subsection{Heegaard Floer homology}

A \textit{Heegaard diagram} for a closed $3$-manifold $M$ is a Heegaard surface $\Sigma$ of genus $g$ together with two systems $\alpha$ and $\beta$ of simple closed curves $\alpha_1$, $\dots$, $\alpha_g$ and $\beta_1$, $\dots$, $\beta_g$ on $\Sigma$ representing generators of $H_1(\Sigma; \Z)$ and bounding disks in the two handlebodies in $M$ respectively.
Heegaard Floer homology is constructed by taking a Heegaard diagram $(\Sigma, \alpha, \beta)$ for $M$ and applying Lagrangian intersection Floer theory~\cite{Fl88, FOOO09} to the tori $\alpha_1 \times \cdots \times \alpha_g$ and $\beta_1 \times \cdots \times \beta_g$ in the symmetric product of $g$ copies of $\Sigma$~\cite{OS04d, OS04e}.

Heegaard Floer homology assigns to $M$ a finitely generated abelian group $\widehat{HF}(M)$ and finitely generated $\Z[U]$-modules $HF^\infty(M)$, $HF^+(M)$, $HF^-(M)$, where $U$ is a formal variable in the polynomial ring $\Z[U]$.
Each group $HF^{\circ}(M)$ has the following decomposition over $\Spin^c(M)$:
\[ HF^\circ(M) = \bigoplus_{\mathfrak{s} \in \Spin^c(M)} HF^\circ(M, \mathfrak{s}). \]
Furthermore, each group $HF^{\circ}(M, \mathfrak{s})$ carries an absolute $\Z / 2 \Z$-grading, and we can take the Euler characteristic $\chi(HF^{\circ}(M, \mathfrak{s}))$ with respect to the grading.
These four flavors of Heegaard Floer homology are related by the following exact triangles:
\begin{align*}
\cdots \to HF^-(M, \mathfrak{s}) \to HF^\infty(M, \mathfrak{s}) \to HF^+(M, \mathfrak{s}) \to HF^-(M, \mathfrak{s}) \to \cdots, \\
\cdots \to \widehat{HF}(M, \mathfrak{s}) \to HF^+(M, \mathfrak{s}) \to HF^+(M, \mathfrak{s}) \to \widehat{HF}(M, \mathfrak{s}) \to \cdots, \\
\cdots \to HF^-(M, \mathfrak{s}) \to HF^-(M, \mathfrak{s}) \to \widehat{HF}(M, \mathfrak{s}) \to HF^-(M, \mathfrak{s}) \to \cdots.
\end{align*}

Ozsv\'ath and Szab\'o~\cite[Theorem 1.2]{OS04e} showed that $HF^+(M)$ is a categolification of Turaev's torsion function $T_M$:
\begin{thm}[\cite{OS04e}] \label{thm:OS-torsion}
Let $M$ be a closed $3$-manifold and $\mathfrak{s} \in \Spin^c(M)$ such that $c_1(\mathfrak{s})$ is not torsion.
Then
\[ \chi(HF^+(Y, \mathfrak{s})) = \pm T_M(\mathfrak{s}). \]
\end{thm}

When $M$ is a rational homology $3$-sphere, $HF^+(M, \mathfrak{s})$ carries an absolute $\Q$-grading.
The \textit{correction term} or \textit{$d$-invariant} $d(M, \mathfrak{s})$ is defined to be the minimal grading of nontorsion elements in the image of the map $\pi \colon HF^\infty(M, \mathfrak{s}) \to HF^+(M, \mathfrak{s})$.
Ozsv\'ath and Szab\'o~\cite[Theorem 1.3]{OS03a} showed that $d(M, \mathfrak{s})$ and $\chi(\coker \pi)$ determine the Casson invariant of an integral homology sphere $M$.

Ozsv\'ath and Szab\'o~\cite[Theorem 1.1]{OS04b} also showed that $\widehat{HF}(M)$ detects the Thurston norm.
\begin{thm}[\cite{OS04b}] \label{thm:OS-norm}
Let $M$ be a closed $3$-manifold $M$.
Then
\[ x_M(\phi) = \min \{ |\langle c_1(\mathfrak{s}) \cup \phi, [M] \rangle| ~;~ \mathfrak{s} \in \Spin^c(M) ~\text{such that}~ \widehat{HF}(M, \mathfrak{s}) \neq 0 \}, \]
for $\phi \in H^1(M; \R)$.
\end{thm}

Ni~\cite[Theorem 1.1]{Ni09a} showed that $HF^+(M)$ detects fiberedness of $M$:
\begin{thm}[\cite{Ni09a}] \label{thm:Ni-fiberedness}
Let $M$ be a closed irreducible $3$-manifold and $S$ a properly embedded surface in $M$ of negative Euler characteristic.
If the group
\[ \bigoplus_{\mathfrak{s} \in \Spin^c(M) ~\text{with}~ \langle c_1(\mathfrak{s}), [S] \rangle = \chi_-(S)} HF^+(M, \mathfrak{s}) \]
is isomorphic to $\Z$, then $M$ fibers over a circle with fiber $S$.
\end{thm}

Theorems~\ref{thm:OS-torsion}, \ref{thm:OS-norm}, \ref{thm:Ni-fiberedness} recover Theorem~\ref{thm:torsion-norm} on the abelian torsion.

\subsection{Knot Floer homology}

Ozsv\'ath and Szab\'o~\cite{OS04c}, and Rasmussen~\cite{Ra03} independently defined the \textit{knot Floer homology} $\widehat{HFK}(L)$ for a null-homologous link $L$ in a closed $3$-manifold $M$.
This finitely generated abelian group $\widehat{HFK}(L)$ refines $\widehat{HF}(M)$ in the sense that there exists a spectral sequence from $\widehat{HFK}(L)$ converging to $\widehat{HF(M)}$.

In the case of a knot $K$ in $S^3$, $\widehat{HFK}(K)$ is bigraded:
\[ \widehat{HFK}(K) = \bigoplus_{(i, j) \in \Z^2} \widehat{HFK}_j(K, i), \]
where $i$ and $j$ are called the \textit{Alexander grading} and the \textit{homological grading} respectively.
The group $\widehat{HFK}(K)$ is a categorification of the classical Alexander polynomial $\Delta_K(t)$ of $K$:
\begin{thm}[\cite{OS04c, Ra03}]
For a knot $K$ in $S^3$ we have
\[ \Delta_K(t) = \sum_{(i, j) \in \Z^2} (-1)^j \left( \rank \widehat{HFK}_j(K, i) \right) t^i. \] 
\end{thm}
As shown in \cite{OS03b}, $\widehat{HFK}(K)$ of an alternating knot $K$ is completely determined by $\Delta_K(t)$ and the signature of $K$.
See also \cite{Ra02}.

As a consequence of the proof of Theorem~\ref{thm:OS-norm}, Ozsv\'ath and Szab\'o~\cite[Theorem 1.2]{OS04b} also showed that $\widehat{HFK}(K)$ determines the knot genus $g(K)$:
\begin{thm}[\cite{OS04b}] \label{thm:HFK1}
For a knot $K$ in $S^3$ we have
\[ g(K) = \max \{ i \in \Z ~;~ \bigoplus_{j \in \Z} \widehat{HFK}_j(K, i) \neq 0 \}. \]
\end{thm}

Theorem~\ref{thm:HFK1} implies that $\widehat{HFK}(K)$ detects the unknot.
See \cite{OS08, Ni09b} for the case of links.
In particular, Ozsv\'ath and Szab\'o~\cite{OS08} showed that the Thurston and Alexander norms agree for the complements of alternating links in $S^3$.

Results by Ghiggini~\cite{Gh08}, Ni~\cite{Ni07}, and Juh\'asz~\cite{Ju08, Ju10} showed that $\widehat{HFK}(K)$ detects fiberedness of $K$:
\begin{thm}[\cite{Gh08, Ni07, Ju08, Ju10}]
A knot $K$ in $S^3$ is fibered if and only if 
\[ \bigoplus_{j \in \Z} \widehat{HFK}_j(K, g(K)) \]
is isomorphic to $\Z$.
\end{thm}

More generally, Juh\'asz~\cite{Ju06} introduced the sutured Floer homology $SFH(M, \gamma)$ for sutured manifolds $(M, \gamma)$ satisfying certain conditions, generalizing $\widehat{HF}(M)$ and $\widehat{HFK}(K)$.
See \cite{FJR11} for the decategorification of $SFH(M, \gamma)$ and the relationship between $SFH(M, \gamma)$ and the Thurston norm for sutured manifolds.
See also \cite{Al14}.

Kronheimer and Mrowka~\cite{KrMr10a, KrMr10b, KrMr11} extended instanton and monopole Floer homology to sutured manifolds.
They also introduced a knot invariant $KHI(K)$ being a categorification of the classical Alexander polynomial and detecting the knot genus and fiberedness of knots.

%%%%%%%% Section 6 %%%%%%%%%%%%%%%%%%%%%%%%%%%%%%%%%%%%%%%%%%%%%%%%%%%%%%%%%%%%
\section{Torsion invariants} \label{sec:6}

The Alexander polynomial has been generalized in different ways to three flavors of nonabelian Alexander polynomials:
twisted Alexander polynomials introduced by Lin~\cite{Linx01} and Wada~\cite{Wa94}, higher-order Alexander polynomials by Cochran~\cite{Co04} and Harvey~\cite{Ha05}, and $L^2$-Alexander invariant by Li and Zhang~\cite{LZ06}.
As seen in the equivalence of the Alexander polynomial and abelian torsion (Theorem~\ref{thm:abelian_torsion}), these generalized Alexander polynomials are also systematically studied in terms of nonabelian Reidemeister torsion:
Reidemeister torsion associated with linear representations, higher-order Reidemeister torsion introduced by Friedl~\cite{Fril07}, and $L^2$-Alexander torsion by Dubois, Friedl and L\"uck~\cite{DFL16}.
We describe relationships between these invariants and the Thurston norm.

\subsection{Reidemeister torsion} \label{ss:R-torsion}

We briefly review Reidemeister torsion associated with linear representations.
See \cite{Mi66, Nic03, Tu01, Tu02b} for details on Reidemeister torsion.

Let $C_* = (C_n \xrightarrow{\partial_n} C_{n-1} \xrightarrow{} \cdots \xrightarrow{} C_0)$ be a finite-dimensional acyclic chain complex over a commutative field $\F$, and let $c = \{ c_i \}$ be a basis of $C_*$.
We choose a basis $b_i$ of $\im \partial_{i+1}$ for each $i$.
Taking a lift $\tilde{b}_{i-1}$ of $\tilde{b}_{i-1}$ in $C_i$ and combining it with $b_i$, we have a basis $b_i \tilde{b}_{i-1}$ of $C_i$ for each $i$.
The \textit{algebraic torsion} $\tau(C_*, c) \in \F \setminus \{ 0 \}$ is defined as:
\[ \tau(C_*, c) = \prod_{i=0}^n [b_i \tilde{b}_{i-1} / c_i]^{(-1)^{i+1}}, \]
where $[b_i \tilde{b}_{i-1} / c_i]$ is the determinant of the base change matrix from $c_i$ to $b_i \tilde{b}_{i-1}$.
It can be checked that $\tau(C_*, c)$ does not depend on the choice of $b_i$ and $b_i \tilde{b}_{i-1}$.

Let $X$ be a connected CW-complex and $R$ a noetherian UFD (e.g., $R$ equals $\Z$ or $\F$).
The cellular chain complex $C_*(\widetilde{X})$ of its universal cover $\widetilde{X}$ is a left $\Z[\pi_1 X]$-module.
We think of $C_*(\widetilde{X})$ also as a right $\Z[\pi_1 X]$-module, using the involution of $\Z[\pi_1 X]$ reversing elements of $\pi_1 X$.
Let $\rho \colon \pi_1 X \to \GL(n, R)$ be a representation.
For each nonnegative integer $i$ the \textit{$i$-th twisted homology group} $H_1^\rho(X; \R^n)$ is defined as:
\[ H_i^\rho(X; R^n) = H_i(C_*(\widetilde{X}) \otimes_{\Z[\pi_1 X]} R^n). \]

The Reidemeister torsion $\tau_\rho(X) \in \F$ associated with a representation $\rho \colon \pi_1 X \to GL(n, \F)$ is defined as follows.
If $H_*^\rho(X; \F^n)$ does not vanish, then we define $\tau_\rho(X) = 0$.
Otherwise, we choose a lift $\tilde{e}$ in $\widetilde{X}$ of each cell $e$ of $X$, and define
\[ \tau_\rho(X) = \tau \left( C_*(\widetilde{X}) \otimes_{\Z[\pi_1 X]} \F^n, \{ \tilde{e} \otimes f_j \}_{e, 1 \leq j \leq n} \right), \]
where $f_1, \dots, f_n$ is the standard basis of $\F^n$.
It is known that $\tau_\rho(X)$ is well-defined as a simple homotopy invariant up to multiplication by elements of $(\pm 1)^n \det \rho(\pi_1 X)$.
Reidemeister torsion $\tau_\rho(X)$ is invariant under conjugation of representations $\rho$.

\begin{rem}
Turaev introduced a refinement of Reidemeister torsion $\tau_\rho(X)$ as an element of $\F$ without any indeterminacy, by fixing an orientation of $H_*(X; \R)$ and an Euler structure of $X$, which is an equivalence class of the choice of lifts $\tilde{e}$~\cite{Tu01, Tu02b}.
\end{rem}

\subsection{Twisted Alexander polynomials}

We describe the results by Friedl and Vidussi~\cite{FV08, FV11a, FV11c, FV14a, FV15}, and Friedl and Nagel~\cite{FN15} that twisted Alexander polynomials detect the Thurston norm and fiberedness of a $3$-manifold.
For more details on twisted Alexander polynomials we refer the reader to the survey papers \cite{DFL15, FV11b, Mori15}.

Let $M$ be a $3$-manifold with empty or toroidal boundary, $\psi \colon \pi_1 M \to F$ a homomorphism to a free abelian group $F$ and $\rho \colon \pi_1 X \to \GL(n, R)$ a representation over a noetherian UFD $R$.
We write $\psi \otimes \rho \colon \pi_1 M \to \GL(n, R[F])$ for the tensor representation defined by $\psi \otimes \rho(\gamma) = \psi(\gamma) \rho(\gamma)$ for $\gamma \in \pi_1 M$.
Then $H_i^{\psi \otimes \rho}(M; R[F]^n)$ is a finitely generated $R[F]$-module for each $i$.
The \textit{$i$-th twisted Alexander polynomial} $\Delta_{M, i}^{\psi, \rho} \in R[F]$ of $M$ associated with $\psi$ and $\rho$ is defined to be its order over $R[F]$, which is well-defined up to multiplication by units in $R[F]$.
We set $\Delta_M^{\psi, \rho} = \Delta_{M, 1}^{\psi, \rho}$.
Twisted Alexander polynomials $\Delta_{M, i}^{\psi, \rho}$ are invariant under conjugation of representations $\rho$.

When $\psi \colon \pi_1 M \to H_1(M)_f$ is the canonical projection and $\rho$ is the trivial representation, $\Delta_M^{\psi, \rho}$ coincides with the usual Alexander polynomial $\Delta_M$ of $M$ defined in Section~\ref{ss:Alexander_polynomial}.
We identify $H^1(M; \Z)$ with $\hom(\pi_1 M, \Z)$ and $R[\Z]$ with the polynomial ring $R[t, t^{-1}]$ so that $1 \in \Z$ corresponds to $t$.
Then for $\phi \in H^1(M; \Z)$ twisted Alexander polynomials $\Delta_{M, i}^{\phi, \rho}$ are in $R[t, t^{-1}]$.

The Reidemeister torsion $\tau_{\psi \otimes \rho}(M) \in Q(R)(F)$ associated with $\psi \otimes \rho \colon \pi_1 M \to GL(n, Q(R)(F))$ is defined as in Section~\ref{ss:R-torsion}, where $Q(R)(F)$ is the quotient field of $R[F]$.
For $\phi \in H^1(M; \Z)$, $\tau_{\phi \otimes \rho}(M) \in Q(R)(t)$ is also defined.
We define
\[\deg \left( a_l t^l + a_{l+1} t^{l+1} + \cdots + a_m t^m \right) = m - l \]
for $l$, $m \in \Z$ with $l < m$ and $a_i \in R$ with $a_l a_m \neq 0$, and further define
\[ \deg \frac{p(t)}{q(t)} = \deg p(t) - \deg q(t) \]
for $p(t), q(t) \in R[t, t^{-1}] \setminus \{ 0 \}$.

The following is a relationship between Reidemeister torsion and twisted Alexander polynomials.
See also \cite{KL99, Kito96}.

\begin{prop}[\cite{FrKim06, FrKim08b, Tu01}] \label{prop:torsion}
Let $M$ be a $3$-manifold with empty or toroidal boundary.
For a homomorphism $\psi \colon \pi_1 M \to F$ and a representation $\rho \colon \pi_1 M \to \GL(n, R)$, if $\Delta_M^{\psi, \rho} \neq 0$, then
\[ \tau_{\psi \otimes \rho}(M) = \frac{\Delta_{M}^{\psi, \rho}}{\Delta_{M, 0}^{\psi, \rho} \Delta_{M, 2}^{\psi, \rho}}. \]
\end{prop}

\begin{rem}
It can be checked for any $\psi$ and $\rho$ that $\Delta_{M, 0}^{\psi, \rho} \neq 0$ and $\Delta_{M, 3}^{\psi, \rho} = 1$, and that $\Delta_{M, 2}^{\psi, \rho} = 0$ if and only if $\Delta_M^{\psi, \rho} = 0$.
The second one follows from the first one and an Euler characteristic argument.
See \cite[Proposition 3.2]{FV11b}.
\end{rem}

Proposition~\ref{prop:torsion}, in particular, shows that $\Delta_M^{\psi, \rho} = 0$ if and only if $H_*^{\psi \otimes \rho}(M; Q(R)(F)^n) = 0$, and that for $\phi \in H^1(M; \Z)$, if $\Delta_M^{\phi, \rho} \neq 0$, then
\[ \deg \tau_{\phi \otimes \rho}(M) = \deg \Delta_M^{\phi, \rho} - \deg \Delta_{M, 0}^{\phi, \rho} - \deg \Delta_{M, 2}^{\phi, \rho}. \]

An advantage of twisted Alexander polynomials and the corresponding Reidemeister torsion is that if a representation $\rho$ is given explicitly, then these invariants can be combinatorially computed, for example, by Fox derivatives for a presentation of the fundamental group.

Friedl and Kim~\cite[Theorems 1.1, 1.2]{FrKim06} generalized Theorem~\ref{thm:torsion-norm} to twisted Alexander polynomials:
\begin{thm}[\cite{FrKim06}] \label{thm:TAP1}
Let $M$ be a $3$-manifold with empty or toroidal boundary and $\rho \colon \pi_1 M \to \GL(n, R)$ a representation.
For $\phi \in H^1(M; \Z)$, if $\Delta_M^{\phi, \rho} \neq 0$, then
\[ x_M(M) \geq \frac{1}{n} \deg \tau_{\phi \otimes \rho}(M). \]
Furthermore, $\Delta_M^{\phi, \rho} \neq 0$ and equality holds if $\phi$ is fibered with $M \neq S^1 \times D^2$ and $M \neq S^1 \times D^2$.
\end{thm}

\begin{rem} \label{rem:monic}
Friedl and Kim~\cite[Theorems 1.3]{FrKim06} also showed that under the assumptions of Theorem~\ref{thm:TAP1}, if $\phi \in H^1(M; \Z)$ is fibered, then $\Delta_M^{\phi, \rho}$ is monic, i.e., its highest and lowest coefficients are units in $R$.
As explained in \cite[Proposition 6.1]{FrKim06}, this theorem can be deduced from Theorem~\ref{thm:TAP1}.
\end{rem}

Generalizing the Alexander norm, Friedl and Kim~\cite[Theorems 3.1, 3.2]{FrKim08b} also defined the \textit{twisted Alexander norm} $||\cdot||_A^\rho$ on $H^1(M; \R)$ associated with a representation $\rho \colon \pi_1 M \to \GL(n, R)$ and generalized Theorem~\ref{thm:A-norm}.
Let $\psi \colon \pi_1 M \to H_1(M)_f$ be the canonical projection, and write $\Delta_M^{\psi, \rho} = \sum_{h \in H_1(M)_f} a_h h$ for $a_h \in R$.
If $\Delta_M^{\psi, \rho} = 0$, we define $||\cdot||_A^\rho = 0$.
Otherwise, we define
\[ ||\phi||_A^\rho = \max \{ \langle \phi, h - h' \rangle \mid h, h' \in H_1(M)_f ~\text{such that}~ a_h a_{h'} \neq 0 \} \]
for $\phi \in H^1(M; \R)$.  
It is clear that $||\cdot||_A^\rho$ is a seminorm on $H^1(M; \R)$ for any $\rho$.
When $\rho$ is the trivial representation, $||\cdot||_A^\rho = ||\cdot||_A$.

\begin{thm}[\cite{FrKim08b}] \label{thm:TAP2}
Let $M$ be a $3$-manifold with empty or toroidal boundary and $\rho \colon \pi_1 M \to \GL(n, R)$ a representation.
Suppose that $b_1(M) > 0$.
Then 
\[ x_M(\phi) \geq \frac{1}{n} ||\phi||_A^\rho \]
for $\phi \in H^1(M; \R)$.
Furthermore, equality holds for $\phi$ in the cone on a fibered face of $B_M$ with $M \neq S^1 \times S^2$ and $M \neq S^1 \times D^2$.
\end{thm}

Theorems~\ref{thm:TAP1}, \ref{thm:TAP2} also provide a fibering obstruction. 
Fibering obstructions on twisted Alexander polynomials in various level of generality were proved in \cite{Ch03, Fril14, FrKim08b, KiMo05, Kita15, GKM05, Pa06}.

As a corollary of the duality of (refined) Reidemeister torsion, Friedl, Kim and the author~\cite[Theorem 1.4]{FrKiKi12} proved the following theorem.
\begin{thm}[\cite{FrKiKi12}]
Let $M$ be an irreducible $3$-manifold with empty or toroidal boundary such that $M \neq S^1 \times D^2$, and let $\rho \colon \pi_1 M \to U(n)$ be a representation.
For $\phi \in H^1(M; \Z)$ whose restriction to any component of $\partial M$ is nontrivial, if $\Delta_M^{\phi, \rho} \neq 0$, then
\[ \deg \tau_{\phi \otimes \rho}(M) \equiv n x_M(\phi) \mod 2. \]
\end{thm}

Using the virtually special theorem by Agol~\cite{Ag13}, Liu~\cite{Liu13} and Przytycki and Wise~\cite{PW18, Wis21}, Friedl and Vidussi~\cite[Theorem 1.2, Corollary 5.10]{FV15} with an extension by Friedl and Nagel~\cite[Theorem 1.3]{FN15} showed that twisted Alexander polynomials determine the Thurston norm:

\begin{thm}[\cite{FV15, FN15}] \label{thm:FV-norm}
Let $M$ be an irreducible $3$-manifold with empty or toroidal boundary.
Then there exists a representation $\rho \colon \pi_1 M \to \GL(n, \C)$ with finite image such that $\Delta_M^{\phi, \rho} \neq 0$ and
\[ x_M(\phi) = \frac{1}{n} \deg \tau_{\phi \otimes \rho}(M)  \]
for all $\phi \in H^1(M; \Z)$.
\end{thm}

\begin{cor}[\cite{FV15}] \label{cor:FV-norm}
Let $M$ be an irreducible $3$-manifold with empty or toroidal boundary.
Suppose that $b_1(M) > 1$.
Then there exists a representation $\rho \colon \pi_1 M \to \GL(n, \C)$ with finite image such that
\[ x_M(\phi) = \frac{1}{n} ||\phi||_A^\rho \]
for all $\phi \in H^1(M; \R)$.
\end{cor}

As explained in \cite[Section 6]{FV15} Theorem~\ref{thm:FV-norm} gives an effective algorithm to compute the Thurston norm.
We will see another algorithm in terms of normal surface theory in Section~\ref{ss:normal_surface}. 
Also, Theorem~\ref{thm:FV-norm} and Corollary~\ref{cor:FV-norm}, in particular, show that the Thurston norm is an invariant of fundamental groups of $3$-manifolds.
We will discuss more on this point of view in Section~\ref{sec:8}. 

Friedl and Vidussi~\cite{FV08, FV11a, FV11c, FV14a} also showed that twisted Alexander polynomials detect fiberedness of $3$-manifolds.
Based on different ideas using Novikov--Sikorav homology, Sikorav~\cite{Si21} showed the fibering detection theorem for general $\phi \in H^1(M; \R)$.

\begin{thm}[\cite{FV14a}] \label{thm:FV-fiberedness}
Let $M$ be a $3$-manifold with empty or toroidal boundary.
If $\phi \in H^1(M; \Z)$ is not fibered, then there exists a representation $\rho \colon \pi_1 M \to \GL(n, \Z)$ with finite image such that $\Delta_{M}^{\phi, \rho} = 0$.
\end{thm}

As a corollary of the fibering detection, together with the study of the Seiberg--Witten invariants of symplectic $4$-manifolds, Friedl and Vidussi~\cite{FV08, FV11b, FV12, FV14a} further showed that a closed $4$-manifold which carries a free circle action admits a symplectic structure if and only if the orbit $3$-manifold is fibered.
The `if' direction generalizes earlier work of Thurston~\cite{Th76}.
See also \cite{Bo09, HW06}.

\subsection{Higher-order Alexander polynomials}

Cochran~\cite{Co04} and Harvey~\cite{Ha05} introduced higher-order Alexander polynomials, analogues of the Alexander polynomial with coefficients in skew fields and showed that their degrees give lower bounds on the Thurston norm.
Following Friedl~\cite{DFL16, Fril07}, we describe results in terms of corresponding higher-order Reidemeister torsion.

Let $\Gamma$ be a torsion-free elementary-amenable group.
By \cite{DLMSY03, KLM88} $\Z[\Gamma]$ is a right (and left) Ore domain, i.e., $\Z[\Gamma]$ embeds in its classical right ring of quotient $\Q(\Gamma) = \Z[\Gamma] \left( \Z[\Gamma] \setminus \{ 0 \} \right)^{-1}$.
The Dieudonn\'e determinant defines a canonical isomorphism $K_1(\Q(\Gamma)) \to \Q(\Gamma)_{ab}^\times$, where $\Q(\Gamma)_{ab}^\times$ is the abelianization of the multiplicative group $\Q(\Gamma) \setminus \{ 0 \}$.

Let $M$ be a $3$-manifold with empty or toroidal boundary.
We define the \textit{higher-order Reidemeister torsion} $\tau_\rho(M) \in \Q(\Gamma)_{ab}^\times \cup \{ 0 \}$ associated with an epimorphism $\rho \colon \pi_1 M \to \Gamma$ onto a torsion-free elementary-amenable group as follows:
If $H_*^\rho(M; \Q(\Gamma)) \neq 0$, then we define $\tau_\rho(M) = 0$.
Otherwise, we take a CW-complex structure of $M$, choose a lift $\tilde{e}$ of each cell $e$, and define
\[ \tau_\rho(M) = \tau \left( C_*(\widetilde{M}) \otimes_{\Z[\pi_1 M]} \Q(\Gamma), \{ \tilde{e} \otimes 1 \}) \right), \]
where $\tau(\cdot) \in \Q(\Gamma)_{ab}^\times$ is the algebraic torsion defined by replacing the usual determinant by the Dieudonn\'e determinant.
The invariant $\tau_\rho(M)$ is well-defined up to multiplication by elements in $\pm \Gamma$.

A pair $(\rho, \phi)$ of an epimorphism $\pi_1 M \to \Gamma$ onto a torsion-free elementary-amenable group and $\phi \in H^1(M; \Z)$ is \textit{admissible} if there exists a homomorphism $\phi_\Gamma \colon \Gamma \to \Z$ such that $\phi_\Gamma \circ \rho \colon \pi_1 M \to \Z$ coincides with $\phi$ under the identification $\hom(\pi_1 M, \Z)$ with $H^1(M; \Z)$.
For an admissible pair $(\rho, \phi)$ we define $\deg_\phi \colon \Q(\Gamma)_{ab}^\times \cup \{ 0 \} \to \Z \cup \{ -\infty \}$ as follows:
Given $p = \sum_{g \in \Gamma} a_g g \in \Z[\Gamma] \setminus \{ 0 \}$, we set
\[ \deg_\phi p = \max \{ \phi_\Gamma(g) - \phi_\Gamma(g') ~;~ a_g a_{g'} \neq 0 \}, \]
and then define
\[ \deg_\phi p q^{-1} = \deg_\phi p - \deg_\phi q \] 
for $p, q \in \Z[\Gamma] \setminus \{ 0 \}$, which induces a homomorphism $\Q(\Gamma)_{ab}^\times \to \Z$.
We extend this to $\deg_\phi 0 = - \infty$.
Now we have an integer-valued invariant $\deg_\phi \tau_\rho(M)$.

\begin{exmp}
Let $M$ be a $3$-manifold with empty or toroidal boundary.
Examples of admissible pairs for $M$ are given by \textit{rational derived series} introduced by Cochran~\cite{Co04} and Harvey~\cite{Ha05}:
We set $\pi_r^{(0)} = \pi_1 M$ and inductively define
\[ \pi_r^{(i)} = \{ \gamma \in \pi_r^{(i-1)} ~;~ \gamma^k \in [\pi_r^{(i-1)}, \pi_r^{(i-1)}] ~\text{for some nonzero $k \in \Z$} \}. \]
Then for any $n$, $\pi_1 M / \pi_r^{(n)}$ is a poly-torsion-free-abelian group, and is, in particular, a torsion-free elementary-amenable group.
We write $\rho_n \colon \pi_1 M \to \pi_1 M / \pi_r^{(n)}$ for the quotient map. 
Then $(\rho_n, \phi)$ is an admissible pair for any $n$ and $\phi \in H^1(M; \Z)$, and we can define $\tau_{\rho_n}(M)$ and $\deg_\phi \tau_{\rho_n}(M)$.
The invariant $\tau_{\rho_n}(M)$ is called the \textit{higher-order Reidemeister torsion of order $n$}.
The invariant $\tau_{\rho_0}(M)$ of order $0$ coincides with the abelian torsion of $M$.
\end{exmp}

Extending the results by Cochran~\cite[Theorem 7.1]{Co04} and Harvey~\cite[Theorem 10.1]{Ha05}, Friedl and Harvey~~\cite[Theorem 1.2]{Fril07}, \cite[Theorem 3.1]{HF07} proved the following theorem, generalizing Theorem~\ref{thm:torsion-norm}:
\begin{thm}[\cite{Fril07, HF07}] \label{thm:higher-order}
Let $M$ be a $3$-manifold with empty or toroidal boundary, and $(\rho, \phi)$ an admissible pair for $M$.
Then
\[ x_M(\phi) \geq \deg_\phi \tau_\rho(M). \]
Furthermore, equality holds if $\phi$ is fibered with $M \neq S^1 \times S^2$ and $M \neq S^1 \times D^2$.
\end{thm}

Theorem~\ref{thm:higher-order} also provides a fibering obstruction.
Friedl~\cite{Fril17} gave more fibering obstructions on higher-order Reidemeister torsion in terms of Novikov-Sikorav homology.

By the duality of higher-order Reidemeister torsion Friedl and Kim~\cite[Theorem 4.4]{FrKim08a} proved the following theorem:
\begin{thm}[\cite{FrKim08a}]
Let $M$ be a closed $3$-manifold or the complement of a link in $S^3$ and $(\rho, \phi)$ an admissible pair for $M$.
Then
\[ \max \{ \deg_\phi \tau_\rho(M), 0 \} \equiv x_M(\phi) \mod 2. \]
\end{thm}

An advantage of higher-order Alexander polynomials or higher-order Reidemeister torsion is that these invariants have monotonicity concerning epimorphisms $\rho$.
Extending the result by Cochran~\cite[Theorem 5.4]{Co04}, Friedl~\cite[Theorem 1.3]{Fril07} and Harvey~\cite[Theorem 2.2, Corollary 2.10]{Ha06} proved the following:
\begin{thm}[\cite{Fril07, Ha06}] \label{thm:triple}
Let $M$ be a $3$-manifold with empty or toroidal boundary, and $(\rho \colon \pi_1 M \to \Gamma, \phi)$ an admissible pair for $M$.
Let $\varphi \colon \Gamma \to \Gamma'$ be an epimorphism such that $(\varphi \circ \rho, \phi)$ is admissible.
Then
\[ \deg_\phi \tau_{\varphi \circ \rho}(M) \leq \deg_\phi \tau_\rho(M). \]
\end{thm}

\begin{cor}[\cite{Co04, Ha06}]
Let $M$ be a $3$-manifold with empty or toroidal boundary and $\phi \in H^1(M; \Z)$.
Then
\[ \deg_\phi \tau_{\rho_{n-1}}(M) \leq \deg_\phi \tau_{\rho_n}(M) \]
for each positive integer $n$.
\end{cor}

\subsection{$L^2$-Alexander torsion}

Li and Zhang~\cite{LZ06} introduced the $L^2$-Alexander invariant, an $L^2$-analogue of higher-order Alexander polynomials.
Dubois, Friedl and L\"uck~\cite{DFL16} introduced the $L^2$-Alexander torsion, the corresponding Reidemeister tosion generalizing the $L^2$-torsion~\cite{Lu02}.
Later, Friedl and L\"uck~\cite{FL17} introduced the universal $L^2$-torsion, which is a further generalization of the invariants.
For more details on these $L^2$-invariants we refer the reader to the survey papers \cite{DFL15, FLT19, Lu21}.
For basic terminology in $L^2$-theory see \cite{Lu02}.

Let $\Gamma$ be a torsion-free group.
We denote by $L^2(\Gamma)$ the Hilbert space of formal sums $\sum_{\gamma \in _\Gamma} a_\gamma \gamma$ for $a_\gamma \in \C$ such that $\sum_{\gamma \in \Gamma} |a_\gamma|^2 < \infty$.
The group von Neumann algebra $\mathcal{N}(\Gamma)$ of $\Gamma$ is defined to be the algebra of bounded $\Gamma$-equivalent operators on $L^2(\Gamma)$.
We denote by $\mathcal{U}(\Gamma)$ the Ore localization of $\mathcal{N}(\Gamma)$ with respect to the multiplicative subset of nonzero divisors.
(It is the algebra of affiliated operators.)
Now we consider the division closure $\mathcal{D}(\Gamma)$ of $\Z[\Gamma]$ in $\mathcal{U}(\Gamma)$, which is the smallest subring $S$ of $\mathcal{U}(\Gamma)$ containing $\Z[\Gamma]$ so that any element of $S$ invertible in $\mathcal{U}(\Gamma)$ is already invertible in $S$.
We say that $\Gamma$ \textit{satisfies the Atiyah conjecture} if given any $m \times n$ matrix $A$ in $\Z[\Gamma]$ the $L^2$-dimesnion of the kernel of the map $L^2(\Gamma)^m \to L^2(\Gamma)^n$ sending $v \to v A$ for $v \in L^2(\Gamma)^m$ is a natural number.
It is an open question whether all torsion-free groups satisfy the Atiyah conjecture.
By \cite{Linn93} $\Gamma$ satisfies the Atiyah conjecture if and only if $\mathcal{D}(\Gamma)$ is a skew-field.

Let $M$ be an aspherical $3$-manifold with empty or toroidal boundary.
It is one of the consequences of the virtually special theorem~\cite{Ag13, Liu13, PW18, Wis21} that if $M$ is not a closed graph manifold, then $\pi_1 M$ satisfies the Atiyah conjecture.
Suppose that $\pi_1 M$ satisfies the Atiyah conjecture.
We define $\rho_u^{(2)}(\widetilde{M}) \in \mathcal{D}(\pi_1 M)_{ab}^\times$ by
\[ \rho_u^{(2)}(\widetilde{M}) = \tau \left( C_*(\widetilde{M}) \otimes_{\Z[\pi_1 M]} \mathcal{D}(\Gamma), \{ \tilde{e} \otimes 1 \}) \right) \in \mathcal{D}(\pi_1 M)_{ab}^\times \]
as in the definition of higher-order Reidemeister torsion.
The invariant $\rho_u^{(2)}(\widetilde{M})$ is well-defined up to multiplication by elements in $\pm \pi_1 M$ and coincides with the one introduced by Friedl and L\"uck in \cite{FL17}, called \textit{the universal $L^2$-torsion} on $M$.

As described in \cite[Section 2.4]{FL17} the Fuglede--Kadison determinant defines a homomorphism $\det_{\mathcal{N}(\pi_1 M)} \colon \mathcal{D}(\pi_1 M)_{ab}^\times \to \R$.
The image $\det_{\mathcal{N}(\pi_1 M)}(\rho_u^{(2)}(\widetilde{M}))$ coincides with the $L^2$-torsion $\rho^{(2)}(\widetilde{M})$ of $M$, which is equal to
\[ - \frac{1}{6 \pi} \sum \vol M_i, \]
where $M_i$ are the hyperbolic pieces in the JSJ decomposition of $M$.
More generally, for $\phi \in H^1(M; \Z)$ and $t \in \R_{>0}$, we can consider the Fuglede--Kadison determinant $\det_{\mathcal{N}(\pi_1 M), t} \mathcal{D}(\pi_1 M)_{ab}^\times \to \R$ twisted by the character on $H_1(M)_f$ sending $h \in H_1(M)_f$ to $t^{\langle \phi, h \rangle}$.
We define a function $\bar\rho^{(2)}(\widetilde{M}; \phi) \colon \R_{>0} \to \R$ by
\[\bar\rho^{(2)}(\widetilde{M}; \phi)(t) = \det{}_{\mathcal{N}(\pi_1 M), t}(\rho_u^{(2)}(\widetilde{M})). \]
The invariant $\bar\rho^{(2)}(\widetilde{M}; \phi)$ is well-defined up to multiplication by functions of the form $t^r$ for some $r \in \R$, and coincides with the \textit{(full) $L^2$-Alexander torsion} of $M$ introduced by Dubois, Friedl and L\"uck~\cite{DFL16}.
In fact, the $L^2$-Alexander torsion itself is defined for any irreducible $3$-manifold $M$ with empty or toroidal boundary.
The following example complements the case of graph manifolds.

\begin{exmp}
Let $M$ be a graph manifold which is not homeomorphic to $S^1 \times S^2$ nor to $S^1 \times D^2$.
Dubois, Friedl and L\"uck~\cite[Theorem 1.1]{DFL16} showed that $\bar\rho^{(2)}(\widetilde{M}; \phi)$ is represented by the function
\[ \bar\rho(t) =
\begin{cases}
1 &\text{if $t \leq 1$}, \\
t^{x_M(\phi)} &\text{if $t \geq 1$}
\end{cases}
\]
for $\phi \in H^1(M; \Z)$.
See also \cite{Her17}.
\end{exmp}

Based on the virtually special theorem~\cite{Ag13, Liu13, PW18, Wis21}, Friedl and L\"uck~\cite[Theorem 0.1]{FL19b}, and Liu~\cite[Theorem 1.2]{Liu17} independently proved the following theorem:
\begin{thm}[\cite{FL19b, Liu17}] \label{thm:L2}
Let $M$ be an irreducible $3$-manifold with empty or toroidal boundary which is not homeomorphic to $S^1 \times D^2$.
Then the limits $\lim_{t \to \infty} \frac{\bar\rho^{(2)}(\widetilde{M}; \phi)}{\ln(t)}$ and $\lim_{t \to 0} \frac{\bar\rho^{(2)}(\widetilde{M}; \phi)}{\ln(t)}$ exist, and
\[ \lim_{t \to \infty} \frac{\bar\rho^{(2)}(\widetilde{M}; \phi)}{\ln(t)} - \lim_{t \to 0} \frac{\bar\rho^{(2)}(\widetilde{M}; \phi)}{\ln(t)}   = -x_M(\phi)  \]
for $\phi \in H^1(M; \Z)$.
\end{thm}

\begin{rem}
In Theorem~\ref{thm:L2} the difference of the limits on the left hand side can be regarded as the `degree' of the function $\bar\rho^{(2)}(\widetilde{M}; \phi)$.
\end{rem}

See \cite[Theorem 0.2]{FL19a} for a related theorem on the $\phi$-twisted $L^2$-Euler characteristic and the Thurston norm.

We now discuss an equivalence class of a pair of convex polytopes in $H_1(M; \R)$ associated with the universal $L^2$-torsion.

The \textit{Minkowski sum} of convex polytopes $P$ and $Q$ in $H_1(M; \R)$ is defined by
\[ P + Q = \{ p + q ~;~ p \in P ~\text{and}~ q \in Q \}. \] 
Two convex polytopes $P$ and $Q$ in $H_1(M; \R)$ are \textit{translation equivalent} if $Q = P + \{ v \}$ for some $v \in H_1(M; \R)$.
We denote by $\mathfrak{P}(M)$ the set of translation equivalence classes of convex polytopes in $H_1(M; \R)$.
The Minkowski sum induces the structure of a commutative monoid on $\mathfrak{P}(M)$.
We denote by $\mathfrak{G}(M)$ the Grothendieck group of $\mathfrak{P}(M)$.
Let $\varphi \colon \pi_1 M \to H_1(M)_f$ be the canonical projection.
Taking a section $H_1(M)_f \to \pi_1 M$, we can identify $\Z[\pi_1 M]$ with $\Z[\ker \varphi][H_1(M)_f]$.
We define a map $\mathcal{P} \colon \Z[\pi_1 M] \to \mathfrak{P}(M)$ as follows:
For $f = \sum_{h \in H_1(M)_f} a_h h \in \Z[\ker \varphi][H_1(M)_f] \setminus \{ 0 \}$ we define $\mathcal{P}(f)$ to be the convex hull of all $h$ with $a_h \neq 0$ in $H_1(M; \R)$.
The map extends as a homomorphism $\mathcal{P} \colon \mathcal{D}(\pi_1 M) \setminus \{ 0 \} \to \mathfrak{G}(M)$, which further induces a homomorphism $\mathcal{P} \colon \mathcal{D}(\pi_1 M)_{ab}^\times \to \mathfrak{G}(M)$.

Friedl and L\"uck~\cite[Theorem 3.35]{FL17} showed that the universal $L^2$-torsion determines the (dual) Thurston norm ball:
\begin{thm}[\cite{FL17}]
Let $M$ be an aspherical $3$-manifold with empty or toroidal boundary such that $\pi_1 M$ satisfies the Atiyah conjecture.
Then 
\[ [B_M^*] = 2 \cdot \mathcal{P}(\rho_u^{(2)}(\widetilde{M})) \in \mathfrak{G}(M). \]
\end{thm}

More generally, the above construction associates an equivalence class of a pair of convex polytopes $\mathcal{P}(\Gamma)$ also to a torsion-free group $\Gamma$ satisfying the Atiyah conjecture and having a finite classifying space $B \Gamma$.
For a group $\Gamma$ admitting a presentation with two generators and one relator, 
reinterpreting results by Friedl, Schreve and Tillmann~\cite{FST17, FT20}, Friedl, L\"uck and Tillmann~\cite[Theorems 3.1, 5.4]{FLT19} described a combinatorial construction of $\mathcal{P}(\Gamma)$ from such a presentation, and showed that $\mathcal{P}(\Gamma)$ determines the Bieri--Neumann--Strebel invariant of $\Gamma$.
See also \cite{PTSN19} for related results.
Funke and Kielak~\cite{FuKie18} studied $\mathcal{P}(\Gamma)$ and its relationships with the Bieri--Neumann--Strebel invariant and higher-order Alexander polynomials for free-by-cyclic groups $\Gamma$.

%%%%%%%% Section 7 %%%%%%%%%%%%%%%%%%%%%%%%%%%%%%%%%%%%%%%%%%%%%%%%%%%%%%%%%%%%
\section{Triangulations} \label{sec:7}

We discuss relationships between triangulations of a $3$-manifold and its Thurston norm.
There are algorithms to compute the Thurston norm ball and its fibered faces from triangulations in normal surface theory.
Also, a $\Z / 2 \Z$-analogue of the Thurston norm is known to give a lower bound on the minimal number of tetrahedra in triangulations.

\subsection{Thurston norm via normal surfaces} \label{ss:normal_surface}

Algorithms to compute the Thurston norm ball in terms of normal surface theory are given by Tollefson and Wang~\cite{TW94, TW96}, and Cooper and Tillmann~\cite{CT09}.
Here we overview a construction of the Thurston norm ball, following \cite{CT09}.

Let $M$ be a closed irreducible $3$-manifold and $\mathcal{T}$ a triangulation of $M$ with $t$ tetrahedra.
Here we mean triangulations to be more general than simplicial triangulations.
We allow triangulations to have simplices with self-identifications on their boundary.  
A triangulation is called \textit{$0$-efficient} if every normal $2$-sphere bounds a $3$-ball contained in a small neighborhood of a vertex. 
Recall that there are $7$ types of normal discs in a tetrahedron $\Delta$:
$4$ triangles around the vertices of $\Delta$ and $3$ quadrilaterals separating the vertices of $\Delta$ into $2$ pairs.
A normal surface is an embedded surface in $M$ whose intersection with each tetrahedron of $\mathcal{T}$ is a collection of disjoint normal discs.
A fundamental fact in normal surface theory is that every incompressible surface in $M$ is isotopic to a normal surface.

If we also take into account of transverse orientations of normal surfaces, there are $2$ equivalence classes for each type of normal discs.
We first consider the linear subspace of the real vector space of dimension $14 t$ with a basis consisting of the equivalence classes of transversely oriented normal discs in the tetrahedra of $\mathcal{T}$.
We denote by $NS^{\nu}(\mathcal{T})$ the linear subspace defined by the so-called matching equations:
for each equivalence class of transversely oriented arcs $\gamma$ in each triangle shared by $2$ tetrahedra $\Delta_\pm$,
\[ t_- + q_- = t_+ + q_+, \]
where $t_\pm$ and $q_\pm$ are the coefficients of the equivalence classes of transversely oriented triangle and quadrilateral in $\Delta_\pm$ respectively containing $\gamma$ in their boundary.
We denote by $NS_+^{\nu}(\mathcal{T})$ the subset of $NS^{\nu}(\mathcal{T})$ consisting of elements whose coefficients are all nonnegative.
An element of $NS_+^{\nu}(\mathcal{T})$ is \textit{admissible} if at most one type of quadrilateral in each tetrahedron is allowed to have nonzero coefficients.

By the construction every admissible integral point of $NS_+^{\nu}(\mathcal{T})$ is represented by a transversely oriented normal surface in $M$, and there are a linear map $\chi^* \colon NS^{\nu}(\mathcal{T}) \to \R$ and a surjective homomorphism $h \colon NS^{\nu}(\mathcal{T}) \to H_2(M; \R)$ corresponding to the Euler characteristic and homology class of a normal surface respectively~\cite[Lemma 3, Proposition 4]{CT09}.
The set $P(\mathcal{T})$ of all elements of $NS_+^{\nu}(\mathcal{T})$ such that the sum of  the coefficients is equal to $1$ is a compact convex polytope in $NS^{\nu}(\mathcal{T})$.
We define $B(\mathcal{T})$ to be the convex hull of the points $\frac{v}{|\chi^*(v)|}$, where $v$ is an admissible vertex of $P(\mathcal{T})$ satisfying $\chi^*(v) < 0$.
Now we can state the following theorem~\cite[Theorem 5]{CT09}:

\begin{thm}[\cite{CT09}] \label{thm:normal_surface}
Let $M$ be a closed irreducible atoroidal $3$-manifold with $b_1(M) > 0$, and $\mathcal{T}$ a simplicial or $0$-efficient triangulation.
Then $h(B(\mathcal{T}))$ coincides with $B_M$.
\end{thm}

Together with Haken's algorithm to check whether the complement of an open tubular neighborhood of an embedded surface $S$ is homeomorphic to the product $S \times [0, 1]$~\cite{Ma03}, Theorem~\ref{thm:normal_surface} also gives an algorithm to determine the fibered faces of $B_M$~\cite[Algorithm 6]{CT09}. 
An alternative algorithm to construct $B_M$ is given in \cite[Algorithm 5.9]{TW96}.

\subsection{$\Z / 2 \Z$-Thurston norm and complexity of $3$-manifolds} \label{ss:z2}

Jaco, Rubinstein and Tillmann~\cite{JRT13} introduced a $\Z/ 2 \Z$-analogue of the Thurston norm.
Let $M$ be a closed irreducible $3$-manifold.
Every cohomology class in $H^1(M; \Z / 2 \Z)$ is the Poincar\'e dual of the homology class represented by a possibly nonorientable embedded surface with some components in $M$.
The \textit{$\Z / 2 \Z$-Thurston norm} $x_{M, 2}$ is the function on $H^1(M; \Z / 2 \Z)$ defined by
\[ x_{M, 2}(\phi) = \min \{ \chi_-(S) ~;~ \text{$S$ is a possibly nonorientable embedded surface dual to $\phi$} \} \]
for $\phi \in H^1(M; \Z / 2 \Z)$.

The \textit{complexity} $c(M)$ of $M$ is the minimal number of tetrahedra in triangulations of $M$.
The number agrees with the one defined by Matveev~\cite{Matv90} unless $M$ is homeomorphic to $S^3$, $\R P^3$ or $L(3, 1)$. 

Generalizing earlier work~\cite{JRT09}, Jaco, Rubinstein, Spreer and Tillmann~\cite[Theorems 1, 2]{JRT13}, \cite[Theorems 1, 3]{JRST20a}, and Nakamura~\cite[Theorems 1.1, 1.2]{Nak17} showed that $x_{M, 2}$ gives lower bounds on $c(M)$:
\begin{thm}[\cite{JRST20a, Nak17}] \label{thm:JRST1}
Let $M$ be a closed irreducible $3$-manifold not homeomorphic to $\R P^3$.
Then
\[ c(M) \geq 1 + 2 x_{M, 2}(\phi). \]
Furthermore, if equality holds, then $M$ is a lens space. 
\end{thm}

\begin{thm}[\cite{JRT13, JRST20a, Nak17}] \label{thm:JRST2}
Let $M$ be a closed irreducible $3$-manifold and suppose $H^1(M; \Z / 2 \Z)$ contains a subgroup $H$ of rank $2$.
Then
\[ c(M) \geq 2 + \sum_{\phi \in H \setminus \{ 0 \}} x_{M, 2}(\phi). \]
Furthermore, if equality holds, then $M$ is a generalized quaternionic space.
\end{thm}

Jaco, Rubinstein, Spreer and Tillmann~\cite[Theorem 1]{JRST20b} also showed that another $\Z / 2 \Z$-analogue of the Thurston norm gives a lower bound on the minimal number of tetrahedra in ideal triangulations of cusped hyperbolic $3$-manifolds.

For rational homology $3$-spheres $M$, Ni and Wu~\cite[Corollary 1.2]{NW15} gave a lower bound on $x_{M, 2}$ by the $d$-invariant $d(M, \mathfrak{s})$, defined via gradings in Heegaard Floer homology~\cite{OS03a}.

%%%%%%%% Section 8 %%%%%%%%%%%%%%%%%%%%%%%%%%%%%%%%%%%%%%%%%%%%%%%%%%%%%%%%%%%%
\section{Profinite rigidity} \label{sec:8}

What properties of $3$-manifolds are determined by the set of finite quotients of their fundamental groups, or the profinite completions of their fundamental groups?
We describe results on the Thuston norm and fiberedness by Boileau and Friedl~\cite{BoFr19}, Bridson, Reid and Wilton~\cite{BR20, BRW17}, and Liu~\cite{Liu20}.
We refer the reader to the survey \cite{Re18} for recent work on profinite rigidity of residually finite groups. 

The \textit{profinite completion} $\widehat{\pi}$ of a group $\pi$ is defined to be the limit $\lim_{\leftarrow} \pi / \Gamma$ of the inverse system $\{ \pi / \Gamma \}_\Gamma$, where $\Gamma$ runs over all finite index normal subgroups of $\pi$.
It easily follows from consideration of finite abelian quotients that if the profinite completions $\widehat{\pi}$ and $\widehat{\pi}'$ of finitely generated groups $\pi$ and $\pi'$ are isomorphic, then so are $H_1(\pi; \Z)$ and $H_1(\pi'; \Z)$.

Bridson, Reid and Wilton~\cite[Theorem A, Corollary 1.1]{BR20}, \cite[Theorem C]{BRW17} showed the following rigidity theorems on fiberedness of $3$-manifolds $M$ with $b_1(M) = 1$:
\begin{thm}[\cite{BR20}] \label{thm:BR}
Let $M_1$ and $M_2$ be $3$-manifolds with $b_1(M_1) = b_2(M_2) = 1$.
Suppose that $\widehat{\pi_1 M_1}$ and $\widehat{\pi_1 M_2}$ are isomorphic.
If $M_1$ has nonempty incompressible boundary and fibers over a circle such that $\pi_1 M_1$ is isomorphic to a semidirect product of the free group of rank $r$ and $\Z$, then so does $M_2$.  
\end{thm}

\begin{thm}[\cite{BRW17}] \label{thm:BRW}
Let $M_1$ and $M_2$ be $3$-manifolds with $b_1(M_1) = b_2(M_2) = 1$.
Suppose that $\widehat{\pi_1 M_1}$ and $\widehat{\pi_1 M_2}$ are isomorphic.
If $M_1$ is a closed hyperbolic $3$-manifold fibering over a circle with fiber of genus $g$, then so is $M_2$. 
\end{thm}

Using different methods, Boileau and Friedl~\cite[Theorems 1.1, 4.6]{BoFr20} showed the following theorems:
\begin{thm}[\cite{BoFr20}] \label{thm:BF1}
Let $M_1$ and $M_2$ be aspherical $3$-manifolds with empty or toroidal boundary.
Suppose that there exists an isomorphism $\widehat{\pi_1 M_1} \to \widehat{\pi_1 M_2}$ such that the induced isomorphism $\widehat{H_1(M_1; \Z)} \to \widehat{H_1(M_2; \Z)}$ is induced by an isomorphism $f \colon H_1(M_1; \Z) \to H_1(M_2; \Z)$.
Then
\[ x_{M_1}(f^* \phi) = x_{M_2}(\phi) \]
for $\phi \in H_1(M_2; \Z)$.
Furthermore, $f^* \phi$ is fibered if and only if so is $\phi$.
\end{thm}

\begin{thm}[\cite{BoFr20}] \label{thm:BF2}
Let $M_1$ and $M_2$ be aspherical $3$-manifolds with empty or toroidal boundary such that $H_1(M_1; \Z)$ and $H_2(M_2; \Z)$ are infinite cyclic groups.
Let $\phi_1 \in H_1(M_1; \Z)$ and $\phi_2 \in H_1(M_2; \Z)$ be generators.
Suppose that there exists an isomorphism $\widehat{\pi_1 M_1} \to \widehat{\pi_1 M_2}$.
Then
\[ x_{M_1}(\phi_1) = x_{M_2}(\phi_2). \]
Furthermore, $\phi_1$ is fibered if and only if so is $\phi_2$.
\end{thm}

The proofs of Theorems~\ref{thm:BF1}, \ref{thm:BF2} rest on the facts that $3$-manifold groups are good in the sense of Serre~\cite{Ser97}, and that the profinite completion of a $3$-manifold group contains enough information on certain twisted Alexander polynomials determining the Thurston norm and fiberedness as in Theorems~\ref{thm:FV-norm}, \ref{thm:FV-fiberedness}.

For knots in $S^3$ Theorems~\ref{thm:BR}, \ref{thm:BF2} can be restated as follows:
\begin{cor}[\cite{BR20, BoFr20}] \label{cor:BF}
Let $J$ and $K$ be two knots in $S^3$ such that $\widehat{\pi_1 X_J}$ and $\widehat{\pi_1 X_K}$ are isomorphic.
Then $g(J) = g(K)$.
Furthermore, $J$ is fibered if and only if $K$ is fibered.
\end{cor}

Bridson and Reid~\cite[Theorem A, Proposition 3.10]{BR20} showed that the profinite completions of fundamental groups distinguish each of the complements of the trefoil knot and the figure-eight knot, and the Gieseking manifold among $3$-manifolds.
Also, Boileau and Friedl~\cite[Theorem 1.5]{BoFr20} showed that every torus knot is also distinguished among knots in $S^3$.

For examples, Hempel's pairs~\cite{Hem14} give examples of Seifert fibered spaces the profinite completions of whose fundamental groups are isomorphic but not satisfying the assumption as in Theorem~\ref{thm:BF1}.
For hyperbolic $3$-manifolds, Liu~\cite[Theorems 1.2, 1.3]{Liu20} strengthen Theorem~\ref{thm:BF1}:
\begin{thm}[\cite{Liu20}] \label{thm:Liu}
Let $M_1$ and $M_2$ be hyperbolic $3$-manifolds. 
Suppose that there exists an isomorphism $\Phi \colon \widehat{\pi_1 M_1} \to \widehat{\pi_1 M_2}$.
Then the following hold:
\begin{enumerate}
\item There exists a unit $\mu \in \widehat{\Z}$ such that $\Phi_* \colon \widehat{H_1(M)_f} \to \widehat{H_1(M_2)_f}$ is induced by an isomorphism $f \colon H_1(M_1)_f \to H_1(M_2)_f$ composed with the multiplication by $\mu$.
\item We have
\[ x_{M_1}(f^* \phi) = x_{M_2}(\phi) \]
for $\phi \in H^1(M_2; \Z)$.
Furthermore, $f^* \phi$ is fibered if and only if so is $\phi$.
\end{enumerate}
\end{thm}

Using Theorem~\ref{thm:Liu} as one of the key ingredients, Liu\cite[Theorem 1.1]{Liu20} proved the following theorem:
\begin{thm}[\cite{Liu20}]
For the fundamental group $\pi$ of a hyperbolic $3$-manifold there exists only finitely many $3$-manifold groups whose profinite completions are isomorphic to $\widehat{\pi}$.
\end{thm}
 
In~\cite{BoFr19} Boileau and Friedl showed that the Thurston norms of all finite covers of an aspherical $3$-manifold determine whether it is a hyperbolic manifold, a graph manifold, or a mixed manifold, i.e., the JSJ decomposition is nontrivial and contains at least one hyperbolic component.
Ueki~\cite{U18} showed that the Alexander polynomial of a knot in $S^3$ is determined by the profinite completion of its knot group.

%%%%%%%% Section 9 %%%%%%%%%%%%%%%%%%%%%%%%%%%%%%%%%%%%%%%%%%%%%%%%%%%%%%%%%%%%

\section{Conjectures and questions} \label{sec:9}

We conclude by collecting some conjectures and questions on the Thurston norm and related topics.

\subsection{Realization problem}

In \cite[Section 4]{Th86} Thurston already gave a large variety of shapes for (dual) Thurston norm balls.
However, the following naive question has been open since the Thurston norm was introduced:
\begin{q}
Which polyhedrons in $\R^n$ are realized as the (dual) Thurston norm balls of $3$-manifolds?
\end{q}

See \cite[Question 6.11]{FLT19} for a restatement of this question and another one in terms of the universal $L^2$-torsion, and see also \cite{PTSN19} for a related result.

\subsection{Complexity functions for circle bundles}

In Theorem~\ref{thm:complexity2} Friedl and Vidussi showed that for all but finitely many circle bundles $N$ over a non-exceptional $3$-manifold $M$ the complexity function $x_N \colon H_2(N; \Z) \to \Z$ is attained by the Thurston norm $x_M$.
As remarked in \cite[Section 1.3]{FV14b} we can ask whether the theorem holds for all circle bundles:
\begin{q}
Let $M$ be a closed irreducible $3$-manifold which is not a closed graph manifold such that $\Delta_M^\phi \neq 0$ for all nontrivial $\phi \in H^1(M; \Z)$, and let $p \colon N \to M$ be an oriented circle bundle. 
Then does the equality 
\[ x_N(\alpha) = |\alpha \cdot \alpha| + x_M(p_* \alpha) \]
hold for any oriented circle bundle $p \colon N \to M$ and any $\alpha \in H_2(N; \Z)$?
\end{q}

\subsection{Twisted Alexander polynomials for hyperbolic knots}

Let $K$ be a hyperbolic knot in $S^3$ and $\phi \in H^1(X_K; \Z)$ a generator.
A holonomy representation $\rho \colon \pi_1 X_K \to \PSL(2, \C)$ of the hyperbolic structure has a lift $\tilde{\rho} \colon \pi_1 X_K \to \SL(2, \C)$~\cite{Cu86, Th80}.
Thus Reidemeister torsion $\tau_{\phi \otimes \tilde{\rho}}(X_K) \in \C(t)$ is defined, and can be checked to be in $\C[t, t^{-1}]$.
Considering Turaev's refinement of $\tau_{\phi \otimes \tilde{\rho}}(X_K)$, Dunfield, Friedl and Jackson~\cite{DFJ12} introduced \textit{hyperbolic torsion polynomial} $\mathcal{T}_K \in \C[t, t^{-1}]$ without any indeterminacy.

Based on experimental results for knots with at most $15$ crossings, Dunfield, Friedl and Jackson~\cite[Conjecture 1.7]{DFJ12} proposed the follwing conjecture:
\begin{conj}[\cite{DFJ12}] \label{conj:DFJ}
Let $K$ be a hyperbolic knot in $S^3$. 
Then
\[ \deg \mathcal{T}_K = 4 g(K) - 2. \]
Furthermore, $K$ is fibered if and only if the leading coefficient of $\mathcal{T}_K$ is equal to $1$.
\end{conj}

Morifuji and Tran~\cite{Mori12, MoTr14, MoTr17} showed that Conjecture~\ref{conj:DFJ} holds for a certain class of $2$-bridge knots.
Later, Agol and Dunfield~\cite{AD20} showed that equality in Conjecture~\ref{conj:DFJ} holds for all libroid hyperbolic knots in $S^3$, including all $2$-bridge knots.
The class of libroid knots is closed under Murasugi sum and contains all special arborescent knots obtained from plumbing oriented bands.
See \cite{MoTr17} for a generalization of Conjecture~\ref{conj:DFJ} for links.

\subsection{Higher-order Alexander polynomials and the knot genus}

Theorems~\ref{thm:higher-order}, \ref{thm:triple} naturally raise the question whether the higher-order Reidemeister torsion determines the Thurston norm.
Dubois, Friedl and L\"uck~\cite[Conjecture 4.4]{DFL16} proposed the following conjecture:
\begin{conj}[\cite{DFL16}] \label{conj:higher-order}
Let $K$ be a knot in $S^3$ and $\phi \in H^1(X_K; \Z)$ a generator.
Then there exists an epimorphism $\rho \colon \pi_1 X_K \to \Gamma$ onto a torsion-free elementary-amenable group such that the pair $(\rho, \phi)$ is admissible and 
\[ \deg_\phi \tau_\rho(X_K) = 2 g(K) - 1. \]
\end{conj}

The following theorem proved by Friedl, Schreve and Tillmann~\cite[Theorem 3]{FST17}, in particular, shows that there are `enough' epimorphisms from knot groups onto torsion-free elementary-amenable groups:
\begin{thm}[\cite{FST17}]
Let $M$ be an irreducible $3$-manifold with empty or toroidal boundary, which is not  a closed graph manifold.
Then $\pi_1 M$ is a residually torsion-free elementary amenable group, i.e., for any nontrivial $\gamma \in \pi_1 M$ there exists an epimorphism $\rho \colon \pi_1 M \to \Gamma$ onto a torsion-free elementary-amenable group such that $\rho(\gamma)$ is nontrivial. 
\end{thm}

\subsection{Lower bounds on complexity of $3$-manifolds}

We have seen in Theorems~\ref{thm:JRST1}, \ref{thm:JRST2} that the $\Z / 2 \Z$-Thurston norm gives lower bounds on the complexity of $3$-manifolds.

The following question was asked by Jaco, Rubinstein and Tillmann in \cite[Section 1]{JRT09}.
\begin{q}[\cite{JRT09}]
Determine an effective bound for the complexity of a closed irreducible $3$-manifold $M$ using a rank $k$ subgroup of $H^1(M; \Z / 2 \Z)$ for $k \geq 3$.
\end{q}

\subsection{Thurston norm balls of finite covers}

As described in~\cite[Proposition 5.4.9]{AFW15}, the following is a consequence of the virtually special theorem~\cite{Ag13, Liu13, PW18, Wis21} and the work of Cooper, Long and Reid~\cite{CLR97}:
\begin{thm} \label{thm:virtual_Haken}
Let $M$ be an irreducible $3$-manifold with empty or toroidal boundary which is not a closed graph manifold.
Then for each positive integer $n$ there exists a finite cover $\widetilde{M}$ of $M$ such that $B_{\widetilde{M}}$ has at least $n$ top-dimensional faces.
\end{thm}

The following is a version of Agol's virtual fibering theorem~\cite{Ag08, Ag13} with a generalization by Kielak~\cite{Kie20b}.

\begin{thm}[\cite{Ag08, Ag13, Kie20b}] \label{thm:virtual_fibering}
Let $M$ be an irreducible $3$-manifold with empty or toroidal boundary which is not a closed graph manifold.
Then there exists a finite covering $p \colon \widetilde{M} \to M$ such that for every nonfibered $\phi \in H^1(M; \R)$, the pullback $p^*(\phi)$ lies in the cone on the boundary of a fibered face of $B_{\widetilde{M}}$. 
\end{thm}

As a corollary of Theorems~\ref{thm:virtual_Haken}, \ref{thm:virtual_fibering} we have the following:
\begin{cor}
Let $M$ be an irreducible $3$-manifold with empty or toroidal boundary which is not a closed graph manifold.
Then for each positive integer $n$ there exists a finite cover $\widetilde{M}$ of $M$ such that $B_{\widetilde{M}}$ has at least $n$ fibered faces.
\end{cor}

See also \cite{DR10, LR08} for related results.

The following questions (also for nonpositively curved graph manifolds in Question~\ref{q:AFW}) were asked by Aschenbrenner, Friedl and Wilton in \cite[Questions 7.5.5, 7.5.6]{AFW15}.

\begin{q}[\cite{AFW15}] \label{q:AFW}
Let $M$ be an irreducible $3$-manifold with empty or toroidal boundary which is not a closed graph manifold.
Does there exist a finite cover $\widetilde{M}$ of $M$ such that all top-dimensional faces of $B_{\widetilde{M}}$ are fibered?
\end{q}

\begin{q}[\cite{AFW15}]
Let $M$ be an irreducible $3$-manifold with empty or toroidal boundary which is not a graph manifold.
Does there exist a finite cover $\widetilde{M}$ of $M$ such that at least one top-dimensional face of $B_{\widetilde{M}}$ is not fibered?
\end{q}

%%%%%%%% References %%%%%%%%%%%%%%%%%%%%%%%%%%%%%%%%%%%%%%%%%%%%%%%%%%%%%%%%%%%

\end{document}